\documentclass[11pt,a4paper]{article}
\usepackage[utf8]{inputenc}
\usepackage{amsmath}
\usepackage{amsfonts}
\usepackage{amssymb}
\usepackage{amsthm}
\usepackage{mathtools}
\usepackage{mathrsfs}
\usepackage{cite}
\usepackage{relsize}
\usepackage{stmaryrd}
\usepackage{hyperref}

\newtheorem{thm}{Theorem}

\newtheorem{defn}[thm]{Definition}

\newcommand{\R}{\mathbb{R}}
\newcommand{\Z}{\mathbb{Z}}
\newcommand{\h}{\mathcal{H}}
\newcommand{\T}{\mathbf{T}}
\newcommand{\Y}{\mathbf{Y}}
\newcommand{\W}{\mathbf{W}}
\newcommand{\sk}{\mathbf{sk}}

\renewcommand{\H}{\mathcal{H}}

\title{Sliding minimal cones in the 3-dimensional half-space}
\author{Edoardo Cavallotto}

\begin{document}
\maketitle
\begin{abstract}
Solving the Plateau problem means to find the surface with minimal area among all
surfaces with a given boundary. Part of the problem actually consists of giving a suitable definition to the notions of ``surface'', ``area'' and ``boundary''. The sliding boundary condition has been introduced by David in order to study the boundary regularity of minimal sets. In order to do that an important step is to know the list of minimal boundary cones, that is to say tangent cones on boundary points of  minimal surfaces. In this paper we focus on cones contained in an half-space and whose boundary can slide along the bounding hyperplane. After giving a classification of one-dimensional minimal cones in the half-plane we provide four new two-dimensional minimal cones in the three-dimensional half space (which cannot be obtained as the Cartesian product of the real line with one of the previous cones).\footnote{The author was supported by a doctoral fellowship within the Marie Skłodowska Curie Initial Training Network MAnET (Metric Analysis for Emergent Technologyes).}

\end{abstract}

\tableofcontents

\section{Introduction}

In order to study boundary regularity of minimal surfaces, in \cite{david2014local} David introduced a new notion of boundary, called \emph{sliding boundary}. Loosely speaking the boundary of a surface subject to this condition is not fixed but is allowed to move in a closed set. A physical example where this condition applies is a soap film contained in a tube: the boundary of the film can move along the inner surface of the tube without leaving it. David proved that sliding minimisers (i.e. minimal surfaces in this new setting) are uniformly rectifiable; moreover he proved that under some mild regularity condition of the boundary, the blow-up limit of a sliding minimal set at a boundary point is a sliding minimal cone with respect to a conical boundary. The assumption of more regularity of the boundary (like $C^1$ or rectifiable) provides flatness (everywhere or almost-everywhere) of its blow-up. Therefore an important step toward understanding the behaviour of sliding minimisers close to the boundary is to list the sliding minimal cones with respect to flat boundaries.

In \cite{fang2016holder} Fang proved that given $\Gamma$ a two-dimensional $C^1$ submanifold of $\R^3$, and a two-dimensional set $E$ which is sliding minimal with respect to $\Gamma$, then $E$ is locally biHölder equivalent to a sliding minimal cone under the assumption that $E$ contains $\Gamma$ and stays on one side of it. In the recent paper \cite{fang2017local} Fang improved this result to $C^{1,\alpha}$-regularity at the boundary of sliding minimal sets. We expect a similar result (at least the Hölder part) to hold with the corresponding analogue of our basic problem when $0\le\alpha\le1$ and follow from our description of minimal cones.

Let $d$ and $n$ be two positive integers such that $1\le d\le n-1$. In the following we will assume our sliding minimisers to have locally finite $d$-dimensional Hausdorff measure, to be contained in the $n$-dimensional half-space $\R^n_+:=\{(x_1,...,x_n):x_n\ge0\}$, and to be subject to the sliding boundary condition given by the bounding hyperplane $\Gamma:=\{(x_1,...,x_n):x_n=0\}$.

The functional we want to minimise is the following weighted Hausdorff measure
\begin{equation}
J_\alpha(E):=\mathcal{H}^d(E\setminus\Gamma)+\alpha\mathcal{H}^d(E\cap\Gamma)
\end{equation}
for $E\subset\mathbb{R}^n_+$. This energy is slightly more general than the d-dimensional Hausdorff measure, and it is is related to functionals appearing in capillarity theory and free boundary problems (see e.g. \cite{finn1974capillarity} \cite{giusti1976boundary} \cite{taylor1977boundary} \cite{dephilippis2015regularity}).

Let us now define precisely the notion of sliding minimiser. Given a set $E\subset\R^n$ with locally finite $d$-dimensional Hausdorff measure we can define the class of sliding competitor as follows.

\begin{defn}[Sliding competitor]\label{sliding_competitor}
Let $\Omega$ be a closed subset of $\R^n$, which may coincide with $\R^n$ itself; and let $\Gamma_i\subset\Omega$, for $0\le i\le I$, be finite family of closed sets. For the sake of notation we will set $\Gamma_0:=\Omega$. We say that $F$ is an admissible competitor to $E$ in $\Omega$, with respect to the sliding boundary domains $\{\Gamma_i\}_{0\le i\le I}$, if there exists a continuous function $\phi:[0,1]\times E\to\mathbb{R}^n$ such that (setting $\phi_t(x)=\phi(t,x)$):
\begin{enumerate}
\item $\phi_0$ is the identity;
\item $F=\phi_1(E)$;
\item $\phi_1$ is Lipschitz;
\item Set $W_t:=\{x\in E:\phi(t,x)\neq x\}$ for $t\in[0,1]$ and $W=\cup_{t\in[0,1]}W_t$, there exists a compact set $K$ such that $\phi(W)\subset K\subset \R^n$.
\item for any $0\le i\le I$, if $x\in E\cap\Gamma_i$ then $\varphi_t(x)\in \Gamma_i$ $\forall t\in[0,1]$.
\end{enumerate}
\end{defn}
Therefore $F$ a sliding competitor for $E$ if it is the image of $E$ under a one-parameter family of compact deformations that take place on a compact set and that moves points along the sliding boundary domains but not away from them. For short we will refer to this class of competitor simply as sliding competitor, and to the deformation satisfying the definition as sliding deformations. In our setting $\Gamma_0=\Omega:=\R^n_+$, $I=1$ and $\Gamma_1:=\Gamma=\partial\Omega$.

\begin{defn}[Sliding minimal set]\label{minimo_scorrevole}
A set $E$ is a $J_\alpha$-minimiser in $\Omega$ with respect to the sliding boundary condition given by $\Gamma$, if, for any sliding competitor $F$, we have that $J_\alpha(E\setminus F)\le J_\alpha(F\setminus E)$.
\end{defn}
Let us remark that both $E$ and $F$ only have locally finite Hausdorff measure, therefore both $J_\alpha(E)$ and $J_\alpha(F)$ may be infinite. However, since admissible competitors are images of compact deformations, it follows that $F$ and $E$ coincide outside a compact set, therefore both $F\setminus E$ and $E\setminus F$ are contained in a compact set and we can compare their energies.

\begin{figure}
\centering
\includegraphics[scale=0.3]{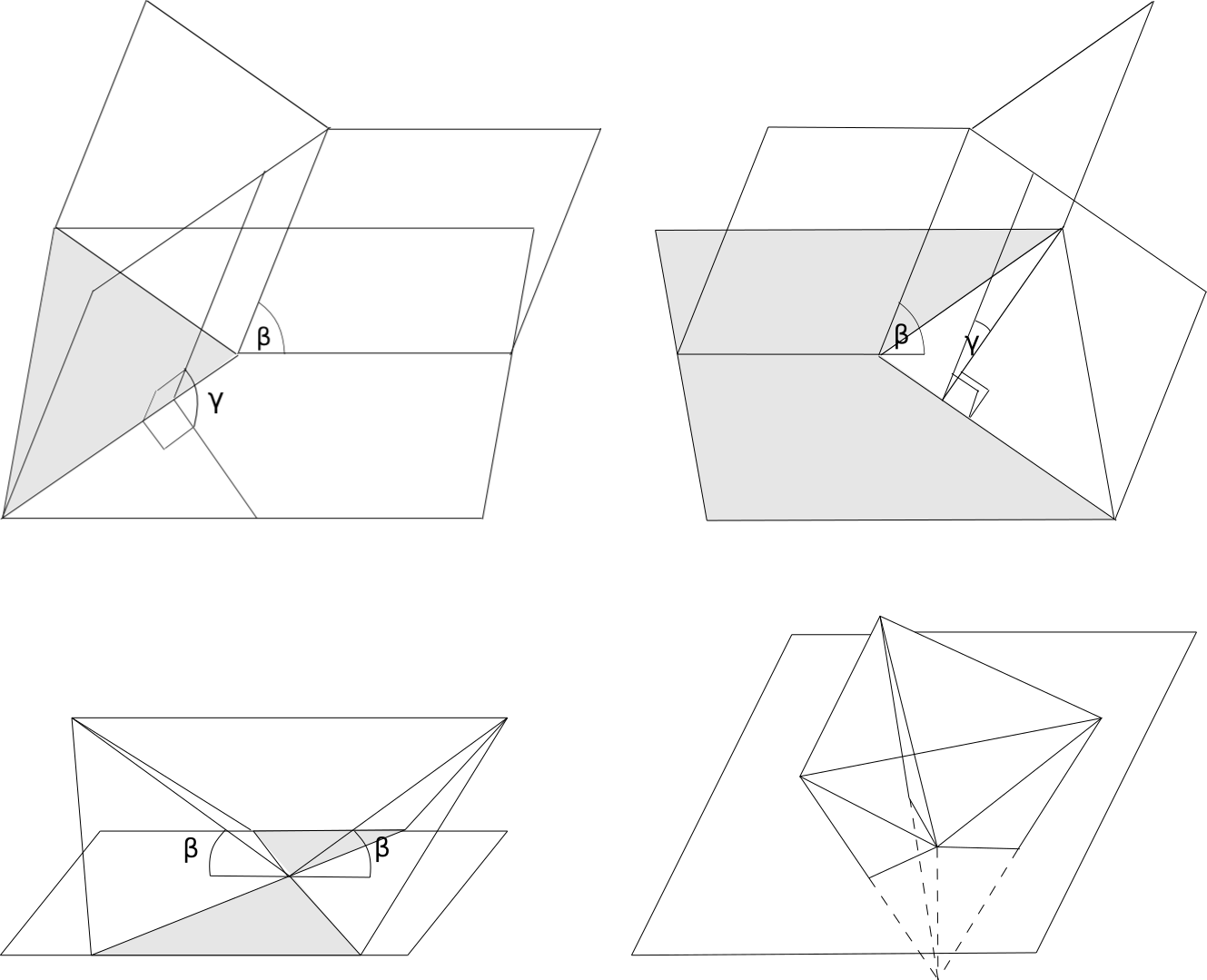}
\caption{The 2-dimensional sliding minimal cones in 3-dimensional half space, respectively: $\Y_\beta$, $\overline{\Y}_\beta$, $\W_\beta$, $\T_+$. For each one of them the grey region represent the intersection of the cone with the horizontal plane.}
\label{coniscorrevoli}
\end{figure}

The purpose of this paper is to study the 2-dimensional minimal cones in the 3-dimensional half-space, that is to say $d=2$ and $n=3$. We will start by classifying the one-dimensional minimal cones in the half-plane, these cones will provide the profiles that a sliding minimal surface can have wile approaching the sliding boundary. Beside the cones obtained as the Cartesian product of $\mathbb{R}$ with a one-dimensional minimal cone in the half-plane, we will prove the sliding minimality of 4 new types of 2-dimensional cones in the 3-dimensional half-space (see Figure \ref{coniscorrevoli}) each of them being indeed a one-parameter family of sliding minimal cones depending on the parameter $\alpha$. In order to prove the sliding minimality of these cones we will use \emph{paired calibrations}. This tool has been employed by Lawlor and Morgan in the proof \cite{lawlor1994paired} of the minimality of the cone over the $(n-2)$-dimensional skeleton of an $n$-dimensional regular simplex; and has recently been generalised to currents with coefficients in a group by Marchese and Massaccesi in \cite{marchese2016steiner} and \cite{massaccesi2014currents}. The technique of paired calibrations consists in applying the divergence theorem to each of the connected components of the complement of a cone using a suitable family of divergence-free vector fields.  

\begin{figure}
\centering
\includegraphics[scale=0.3]{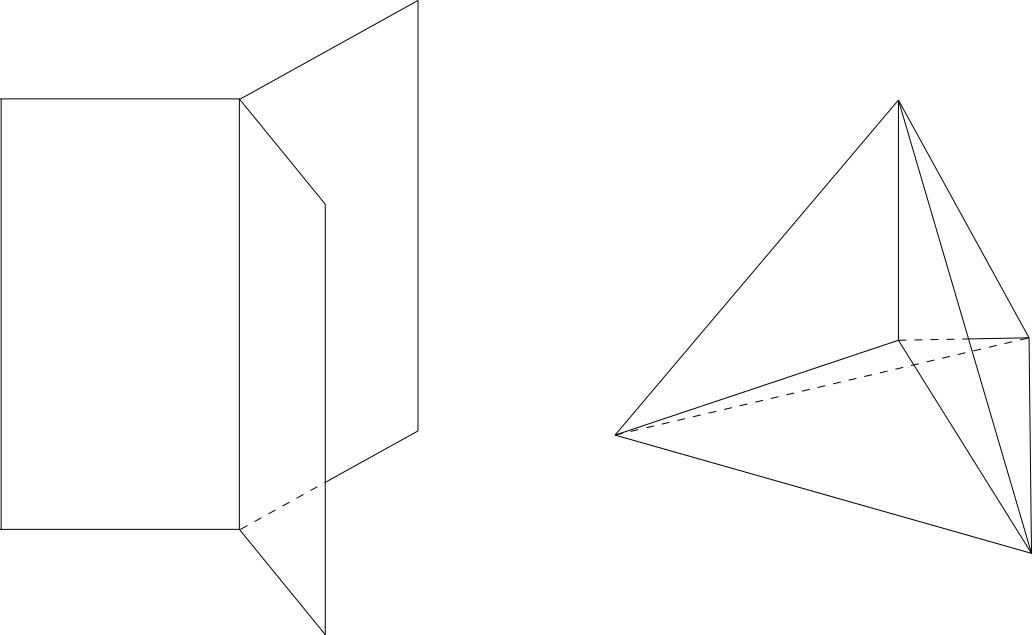}
\caption{Cones of type $\mathbf{Y}$ and $\mathbf{T}$.}
\label{YeTintro.png}
\end{figure}

The complete list of 2-dimensional minimal cones in the whole space $\R^3$ is well known since long time, it was first conjectured by Plateau \cite{plateau1873statique} and then proved by Taylor in \cite{taylor1976structure} a century later (let us also recall the works of Lamarle \cite{lamarle1865stabilite} and Heppes \cite{heppes1964isogonal}). The only minimal cones on the whole space are: planes; the set $\Y$, obtained as the union of three half-panes meeting with equal angle of $120^\circ$; the set $\T$, obtained as the cone over the edges of a regular tetrahedron (see Figure \ref{YeTintro.png}). These cones will be our starting point when looking for sliding minimal cones in the half-space.

The first two examples of minimal cones are obtained with the following procedure. First embed in $\mathbb{R}^3$ a cone of type $\mathbf{Y}$. Then tilt it in such a way that the intersection between the three half-planes meets the horizontal plane with an angle $\beta\in[0,\pi/2]$ and one of the three folds meets orthogonally the horizontal plane. Take now the intersection of this cone with the upper half-space $\mathbb{R}^3_+$. By construction the two sloping folds are forced to meet the horizontal plane with equal angle $\gamma$ (depending on $\beta$), and the vertical fold has the shape of a planar sector whose angle can be either $\beta$ or $\pi-\beta$. In the first case we define $\mathbf{Y}_\beta$ as the union of the cone obtained with the previous construction and the sector of the horizontal plane contained in between the two sloping folds. In the second case we define $\overline{\mathbf{Y}}_\beta$ as the union of the cone obtained with the previous construction and the sector of the horizontal plane not contained in between the two sloping folds. Both $\mathbf{Y}_\beta$ and $\overline{\mathbf{Y}}_\beta$ are minimal if and only if $\cos\gamma=\alpha$.

The third minimal cone is called $\mathbf{W}_\beta$ and is obtained  by taking the union between the intersection of $\overline{\mathbf{Y}}_\beta$ with a half-space bounded by a vertical plane $P$ orthogonal to the vertical fold and its reflection with respect to the plane $P$ itself. The cone $\mathbf{W}_\beta$ is minimal if and only if $\beta\le30^\circ$ and $\cos\gamma=\alpha$.

The fourth minimal cone is called $\mathbf{T}_+$ and is obtained by taking a cone of type $\mathbf{T}$ as in the first picture, flipping it upside down, placing its barycentre at the origin, and finally intersecting it with the half-space $\mathbb{R}^3_+$. Using paired calibrations it is possible to prove the minimality of $\mathbf{T}_+$ for every $\alpha\ge\sqrt{\frac{2}{3}}$. Moreover for every $\alpha<\sqrt{\frac{2}{3}}$ a better competitor to the cone can be found by pinching it down on the horizontal plane in such a way to produce a little triangle and then by connecting it to the boundary in a proper way. 

To know the complete list of sliding minimal boundary cones is the first step in direction of boundary regularity results for sliding minimal sets. We conjecture this list of minimal cones to be complete. However, in order to prove that, it is necessary to classify all the cones satisfying the necessary condition for minimality and, for each one of them, to find a sliding competitor with less energy.

\section{One-dimensional cones in the half-plane}\label{Half-plane}
In this section we will discuss one-dimensional minimal cones in the half-plane. The domain of the sliding boundary will be the bounding line. Using the notation introduced in the previous Section we have that $\Omega=\mathbb{R}^2_+=\{(x,y)\in\mathbb{R}^2:y\ge0\}$ and $\Gamma=\{y=0\}$.

Given $0\le\theta\le\frac{\pi}{2}$, let $P_\theta$ be a half-line meeting $\Gamma$ at the origin with angle $\theta\in[0,\frac{\pi}{2}]$, and let $\theta_\alpha$ be such that $\alpha=\cos\theta_\alpha$. In particular let us remark that $\theta_\alpha\to\pi/2$ when $\alpha\to0$ and $\theta_\alpha\to0$ when $\alpha\to1$. Given $0\le\alpha\le1$ the sliding minimal cones in this setting are the following (see Figure \ref{coni1in2}):

\begin{figure}
\centering
\includegraphics[scale=0.35]{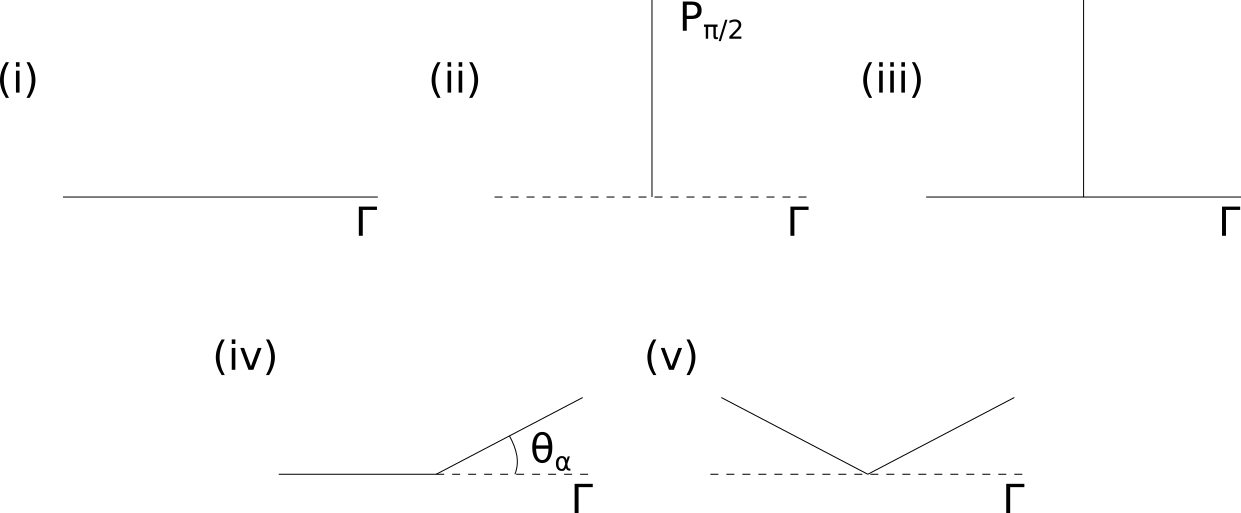}
\caption{One dimensional sliding $\alpha$-minimal cones in the half plane.}
\label{coni1in2}
\end{figure}

\begin{figure}
\centering
\includegraphics[scale=0.35]{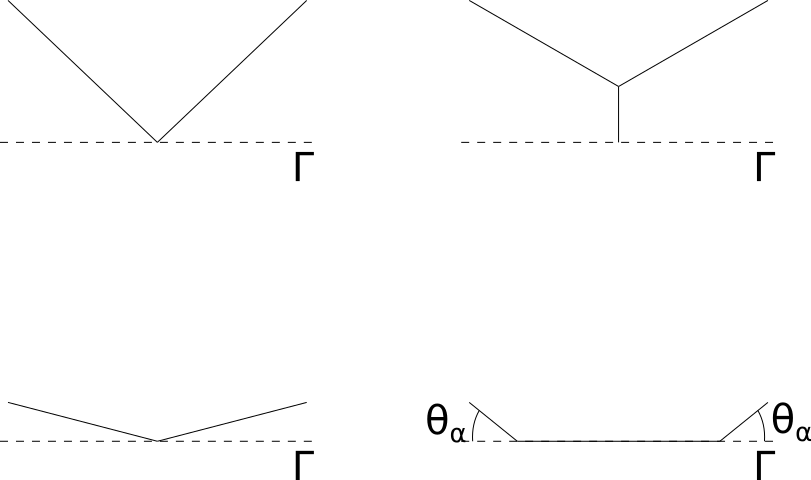}
\caption{Better competitors to the cone (v) in case $\theta>\pi/6$ above, and in case $\theta\le\theta_\alpha$ below.}
\label{competitori_cono_v}
\end{figure}

\begin{itemize}
\item[(i)] $\Gamma$; this cone is trivially minimal because the only member in the family of sliding competitors is $\Gamma$ itself.
\item[(ii)] $P_{\frac{\pi}{2}}$; given a compact set $K$ and $x_0\in P_{\frac{\pi}{2}}\setminus K$, any sliding competitor to the vertical half-line would be (or at least contain) a path connecting $x_0$ to $\Gamma$, therefore it would be longer than the vertical line segment connecting $x_0$ to the origin.
\item[(iii)] $\Gamma\cup P_{\frac{\pi}{2}}$; we can show the minimality of this cone combining the two previous arguments. In fact any competitor to this cone still contains $\Gamma$, and given $x_0$ as before we have that any competitor contains a path connecting $x_0$ to $\Gamma$.
\item[(iv)] the union of $P_{\theta_\alpha}$ with a horizontal half-line laying on $\Gamma$; let $B_1(0)$ be the unit ball centred at the origin, and $\theta\in[0,\frac{\pi}{2}]$. We define $A=(-1,0)$ and $B=(\cos(\theta),\sin(\theta))$ to be the two endpoints of the cone intersected with $B_1(0)$. In order to prove the minimality of the cone it is sufficient to consider all competitors obtained as the union of the two segments $\overline{AC}$ and $\overline{CB}$ where $C=(x,0)$. Hence we have to minimise $J_\alpha$ among a one-parameter family of competitors. Let $E_x$ be one of such competitors, than
\begin{equation}
\begin{aligned}
J_\alpha(E_x) & =\alpha(1+x)+\sqrt{(x-\cos(\theta))^2+\sin^2(\theta)}\\
\left.\frac{\partial J_\alpha(E_x)}{\partial x}\right|_{x=0} & =\left[\alpha+\frac{x-\cos(\theta)}{\sqrt{(x-\cos(\theta))^2+\sin^2(\theta)}}\right]_{x=0}= \alpha-\cos(\theta)\\
\left.\frac{\partial^2 J_\alpha(E_x)}{\partial x^2}\right|_{x=0} & = \left.\frac{\sin^2(\theta)}{\left[(x-\cos(\theta))^2+\sin^2(\theta)\right]^{3/2}}\right|_{x=0}=\sin^2(\theta).
\end{aligned}
\end{equation}
Therefore $x=0$ is a critical point if and only if $\cos(\theta)=\alpha$, and the second derivative is always positive.
\item[(v)] $V_\theta$: the union of $P_\theta$ and its symmetric with respect to the vertical axis, for $\theta_\alpha\le\theta\le\pi/6$. As before let us call $A$ and $B$ the endpoints of $V_\theta$ intersected with the ball $B_1(0)$. Any admissible competitor for $V_\theta$ has to connect $A$, $B$ and $\Gamma$. It is easily seen that $V_\theta$ is minimal among all the competitors obtained as the union of the two segments $\overline{AC}$ and $\overline{CB}$ where $C=(x,0)$. However other types of competitor may occur (see picture \ref{competitori_cono_v}). Pinching together the two segments $\overline{AC}$ and $\overline{CB}$ we can produce a triple junction, and in case $\theta>\pi/6$ it is possible to arrange it in the shape of a $Y$ cone. In this case the obtained competitor is actually the minimiser. Otherwise we can push down the two segments onto $\Gamma$ producing a segment $\overline{C'C''}$ in $\Gamma$. In case $\theta<\theta_\alpha$ we can keep pushing down up to the point when the angles formed by $\overline{AC'}$ and $\overline{C''B}$ with $\Gamma$ turn into $\theta_\alpha$, and again we obtain a minimiser.
\end{itemize}

We remark that cones (i), (ii) and (iii) are independent from $\alpha$ while cone (iv) forms a one parameter family depending on this value. Cones of type (v) form a one-parameter family independent on $\alpha$ but whose endpoints depend on it. When $\alpha=1$ we have that $J_1=\h^1\llcorner\Omega$ and the cone (iv) collapses to the cone (i). On the other side, when $\alpha=0$, we have that $J_0=\h^1\llcorner(\Omega\setminus\Gamma)$. In this case (i) turns into an (even more) trivial minimal cone, and the cones (ii) and (iii) become equivalent with respect to $J_0$ because they only differ from a ``null-measure'' set. Moreover also the cone (iv) collapses to the type (ii)-(iii).

The minimality of the cones (i) and (iv) is only due to the definition of $J_\alpha$ and it would still be minimal without imposing the sliding boundary condition. The cones (ii) and (v) are in the opposite situation: they would not be minimal without the sliding boundary condition, regardless to the coefficient $\alpha$. Cone (iii) in an exceptional case, its minimality relies on the cost functional when $\alpha\le\sqrt{3}/2$ and on the sliding boundary condition when $\alpha>\sqrt{3}/2$. Indeed, in the latter case, what could happen if we drop the sliding boundary condition is that the branching point of the cone could move upwards assuming the $Y$ configuration. Therefore the two branches of this $Y$ would meet the boundary with an angle of $30^\circ$, whose cosine is $\sqrt{3}/2$.

\begin{figure}
\centering
\includegraphics[scale=0.35]{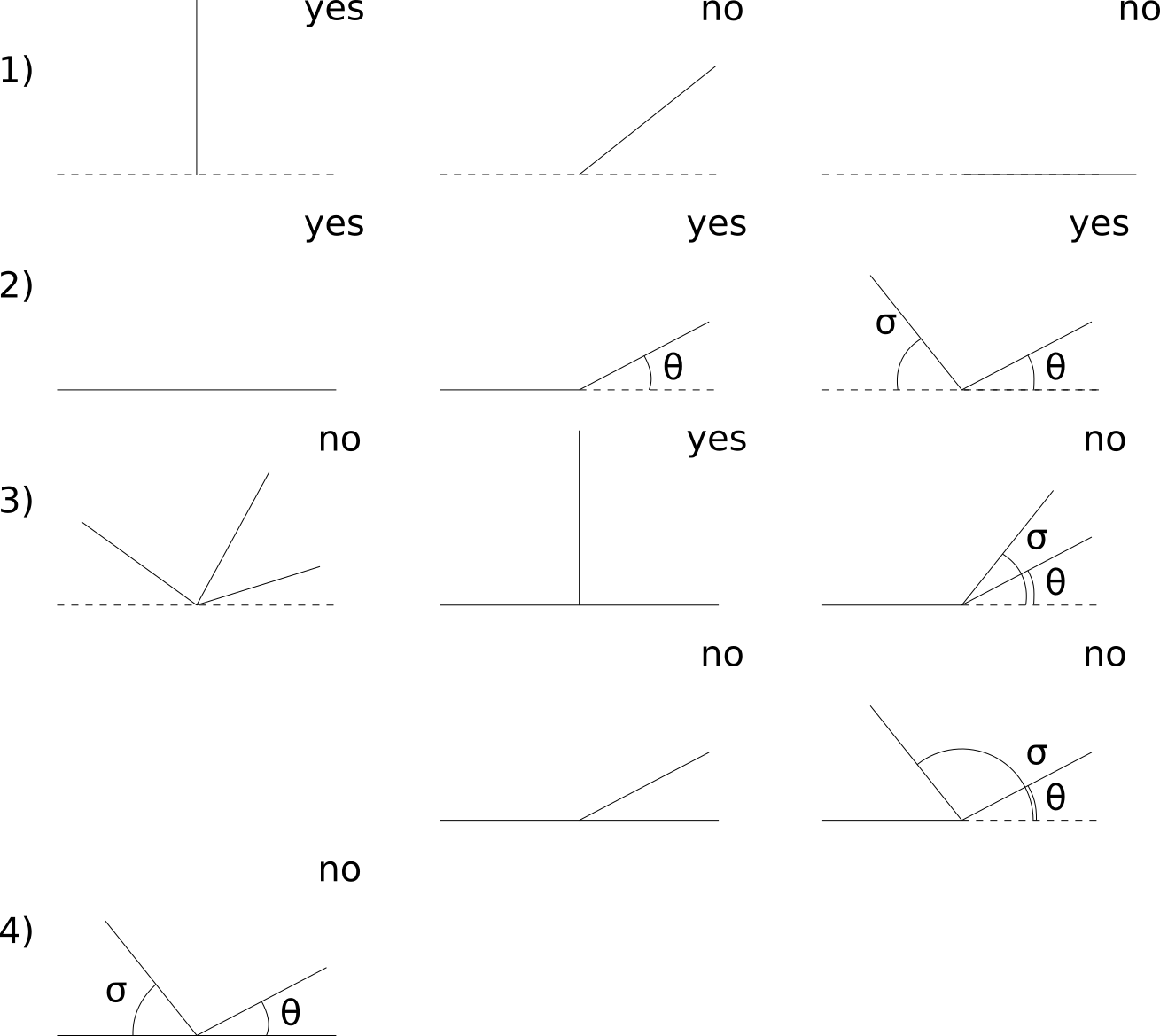}
\caption{Classification of one-dimensional cones by the number of their branches, and whether or not there exists a minimal cone of a given type.}
\label{coni1in2-nonminini}
\end{figure}

In order to show that the aforementioned list of sliding minimal cones is complete we can classify all the one-dimensional cones by the number of distinct half lines (or branches) they are composed by (see Figure \ref{coni1in2-nonminini}) and by their position with respect to $\Gamma$. In case the cone is composed by only one branch than it is sliding minimal only if the branch is vertical (cone of type (ii)); otherwise it is very easy to find a better competitor. In case the cone has two branches than we have three sub-cases, depending on the number of branches contained in $\Gamma$. If they are both contained in $\Gamma$ we find again the cone of type (i); if only one is contained then the cone if minimal if and only if it is of type (iv); and it both the branches are not contained in $\Gamma$ then the cone is minimal if and only if it is of type (v). Let us now discuss the case of three branches. If the three of them are not contained in $\Gamma$ than at least two of them form an angle smaller then $120^\circ$ therefore they can be pinched together decreasing the total length. If only one branch is not contained in $\Gamma$ than the cone is minimal if and only if it is of type (iii) otherwise the sloping branch can be projected onto $\Gamma$ decreasing the total energy. If exactly one branch is contained in $\Gamma$ we can call $\sigma$ and $\theta$ the angles formed by the sloping branches with $\Gamma$ (as in Figure \ref{coni1in2-nonminini}). Than we have two sub-sub-cases. If both the angles are less than a right angle than the angle between them is less than $120^\circ$ and they can be pinched together. If at least one angle is bigger than a right angle than the corresponding branch can be pushed down onto $\Gamma$ in such a way as to obtain a better competitor. The case of four branches is rather simple to rule out. Since there cannot be three branches outside $\Gamma$ (because otherwise we could pinch together two of them as before) the only possibility is that exactly two of them are contained in $\Gamma$, but that means that the other two can be projected onto $\Gamma$ in such a way as to decrease the energy of the set. Therefore there are no four-branched minimal cones. By the same argument no cone with more than 4 branches can be sliding minimal.

\section{Two-dimensional cones in the half-space}\label{two-dimensional cones}

In this Section we will discuss two-dimensional cones in the half-space $\R^3_+:=\{(x,y,z)\in\R^3:z\ge0\}$. The domain of the sliding boundary will be the horizontal plane $\Gamma=\{(x,y,z)\in\R^3:z=0\}$.

\subsection{Cartesian products}\label{Cartesian products}

\begin{figure}
\centering
\includegraphics[scale=0.3]{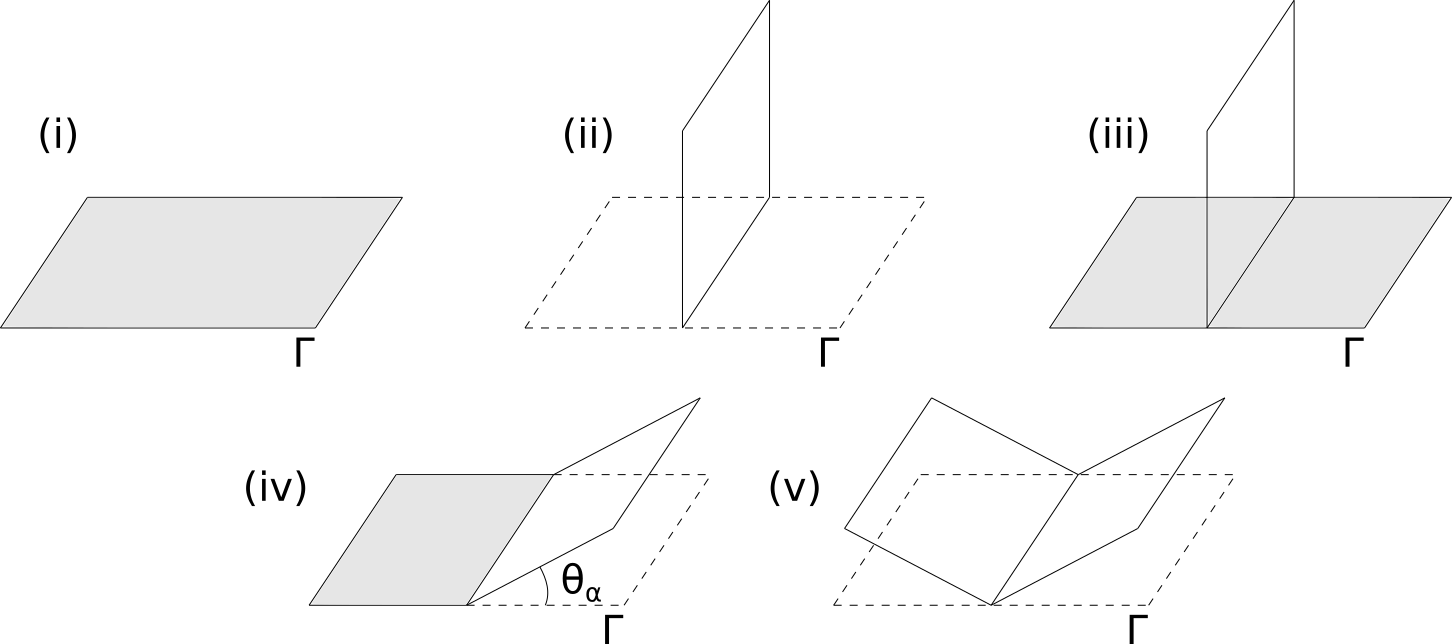}
\caption{Cones obtained as the Cartesian product of $\R$ with a 1-dimensional minimal cone in $\R^2_+$ (the gray region is the intersection between the cones and $\Gamma$).}
\label{coni2in3}
\end{figure}

Let us start our discussion with the 2-dimensional cones that can be obtained as the Cartesian product of $\R$ with one of the 1-dimensional minimal cones in the previous Section (see Figure \ref{coni2in3}).

\begin{figure}
\centering
\includegraphics[scale=0.5]{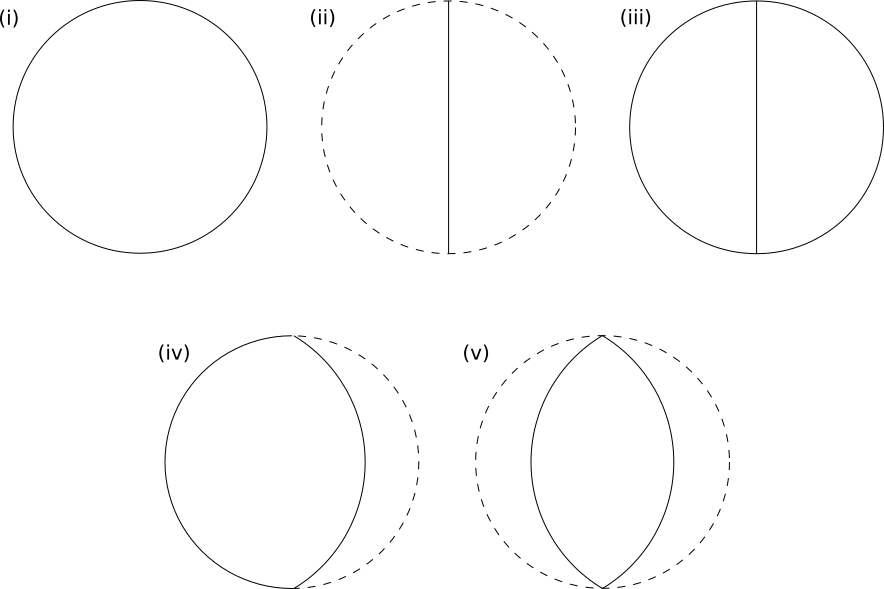}
\caption{Representation of the graphs on the hemisphere (seen from above) characterising the cones in Figure \ref{coni2in3}.}
\label{proditto-cartesiano}
\end{figure}

The minimality of this kind of cones can be proved by a slicing argument. In the following we are going to provide a proof for a cone of type (iv), the minimality of the other cones can be proved in the same way.

Let us denote with $B$ the unit ball of $\R^2$ centred at the origin, and let $A$ be the intersection of the 1-dimensional cone of type (iv) with the ball B. Let $C:=B\times[0,1]\subset \R^2\times\R$ be a cylinder and $D:=A\times[0,1]\subset \R^2\times\R$, then $D$ is the intersection of a 2-dimensional cone of type (iv) with the cylinder $C$. Let us now identify $\R^2\times\R$ with $\R^3$, abusing the notation we will denote with $J_\alpha$ both the functional defined on $\R^2_+$ and the corresponding functional defined on $\R^3_+$. Let $\phi:D\to\R^3$ be a sliding deformation acting in the interior of the cylinder $C$, that is to say $\phi$ is a Lipschitz function and $\phi(W)\subset \textrm{Int}(C)$ where $W:=\{p\in D:\phi(p)\neq p\}$. Therefore $M:=\phi(D)$ is a sliding competitor in the cylinder $C$. Let $\llbracket A\rrbracket\in \mathscr{P}_1(\R^2,\Z_2)$ be the 1-dimensional polyhedral chain with coefficient in $\Z_2$ whose support is $A$, and let $\llbracket D\rrbracket:=\llbracket A\rrbracket\times\llbracket0,1\rrbracket\in \mathscr{P}_2(\R^3,\Z_2)$ be the 2-dimensional polyhedral chain with coefficient in $\Z_2$ whose support is $D$, it follows that $f_\sharp(\llbracket D\rrbracket)$ is supported in $M$ (for the definition and properties of flat chains see, e.g. \cite{federer2014geometric} or \cite{white1999rectifiability}). Let us now define the following orthogonal projection
\begin{equation}
\begin{aligned}
\pi:&\R^2\times\R\to\R\\
&(x,y,t)\mapsto t.
\end{aligned}
\end{equation}
Since $\H^2(M)<+\infty$ we have that $\H^1(M\cap\pi^{-1}(t))<+\infty$ for almost every $t\in[0,1]$, therefore the slice $\langle \phi_\sharp\llbracket D\rrbracket,\pi,t\rangle\in \mathscr{F}_1(\R^2,\Z_2)$ exists for almost every $t\in[0,1]$ since its support is contained in $M\cap\pi^{-1}(t)$.
Moreover the boundary of this slice is the same as the boundary of $\langle \llbracket D\rrbracket,\pi,t\rangle$, and it can be shown using the properties of the slices and the definition of $\phi$ as follows
\begin{equation}
\begin{aligned}
\partial\langle \phi_\sharp\llbracket D\rrbracket,\pi,t\rangle &=\langle \partial\phi_\sharp\llbracket D\rrbracket,\pi,t\rangle\\
&=\langle \phi_\sharp\partial\llbracket D\rrbracket,\pi,t\rangle\\
&=\phi_\sharp\langle\partial\llbracket D\rrbracket,\pi\circ\phi,t\rangle\\
&=\langle\partial\llbracket D\rrbracket,\pi,t\rangle\\
&=\partial\llbracket A\rrbracket.
\end{aligned}
\end{equation}
Hence for almost every $t\in[0,1]$ the slice $M\cap\pi^{-1}(t)$ contains a curve whose length is finite and whose endpoints are the same as the endpoints of $A$, and by the argument in the previous Section we have that $J_\alpha(M\cap\pi^{-1}(t))\ge J_\alpha(A)$. Therefore we can now compute as follows
\begin{equation}
\begin{aligned}
J_\alpha(M)&\ge\int_0^1J_\alpha(M\cap\pi^{-1}(t))d\H^1(t)\ge\int_0^1J_\alpha(A)d\H^1(t)=J_\alpha(D).
\end{aligned}
\end{equation}

\subsection{Characterisation of minimal cones}

Let $S:=\{(x,y,z)\in\R^3:x^2+y^2+z^2=1\}$ be the unit sphere of $\R^3$ centred at the origin, and let $C\subset\R^3_+$ be a sliding minimal cone. We remark that, since a cone is invariant by dilations, it is completely characterised by its intersection with $S$. In our case $C$ is contained in a half-space, hence it is completely characterised by the graph obtained as its intersection with the upper hemisphere $S$, that is to say $S_+:=S\cap\R^3_+$.

In \cite{taylor1973regularity} and \cite{taylor1976structure} Taylor proved that the graph obtained as the intersection of a minimal cone with the unit sphere must consist of arcs of great circles intersecting three at a time at a finite number of points, and the angles of intersection must be $120^\circ$. Moreover she proved that if $A$ is one of the region in which the sphere is divided by the cone, then $A$ is a spherical polygon having at most 5 sides and the lengths of the arcs of these nets can be computed (in terms of the angle at the origin they subtend) using the following formulae:
\begin{itemize}
\item if $A$ is bounded by only one edge then $A$ is a hemisphere bounded by a great circle;
\item if $A$ is bounded by 2 edges then it is a gore whose side length is $\pi$;
\item If $A$ is a spherical equiangular triangle (all angles $120^\circ$), its side
length  is $\arccos(-\frac{1}{3})$;
\item If $A$ is a spherical quadrilateral with $120^\circ$ angles at its vertices, then
it is ``rectangular'' in the sense that opposite sides are of equal length, and
the lengths $\alpha$ and $\beta$ of its adjacent sides are related by the formula
\begin{equation}\label{lati_rettangolo}
\cos(\beta)=\frac{3-5\cos(\alpha)}{5-3\cos(\alpha)}
\end{equation}
which in terms of half angles becomes
\begin{equation}\label{lati_rettangolo/2}
\cos(\beta/2)=2\sin(\alpha/2)\sqrt{1+3\sin^2(\alpha/2)}
\end{equation}
\item If $A$ is a spherical pentagon ($120^\circ$ angles), and if $\alpha$ and $\beta$ are the lengths of adjacent sides, then the length $\gamma$ of the side adjacent to
neither is given by
\begin{equation}\label{lati_pentagono}
2\cos(\gamma)=\frac{1}{3}+\cos(\alpha)+\cos(\beta)+\cos(\alpha)\cos(\beta)-\sin(\alpha)\sin(\beta).
\end{equation}
\end{itemize}
Taylor proved that there are only 10 network satisfying the previous condition on the sphere. She also proved that only 3 of the cones corresponding to these networks are minimal, these cones are: the plane; the cone $\Y$, obtained as the union of three half-planes meeting with equal angle of $120^\circ$; and $\T$, the cone over the edges of a regular tetrahedron (see Figure \ref{YeTintro.png}). Finally Taylor provided a better competitor for each one of the remaining cones.

In our setting we have an extra condition concerning the way in which the network can meet the equator. Applying the slicing argument of the previous subsection to the blow-up of a cone in a point on $\Gamma$, we have that the arcs can meet the equator only with one of the 1-dimensional optimal profiles of Section \ref{Half-plane}.

\subsection{Half $\T$}\label{Half_T}

\begin{figure}
\centering
\includegraphics[scale=0.4]{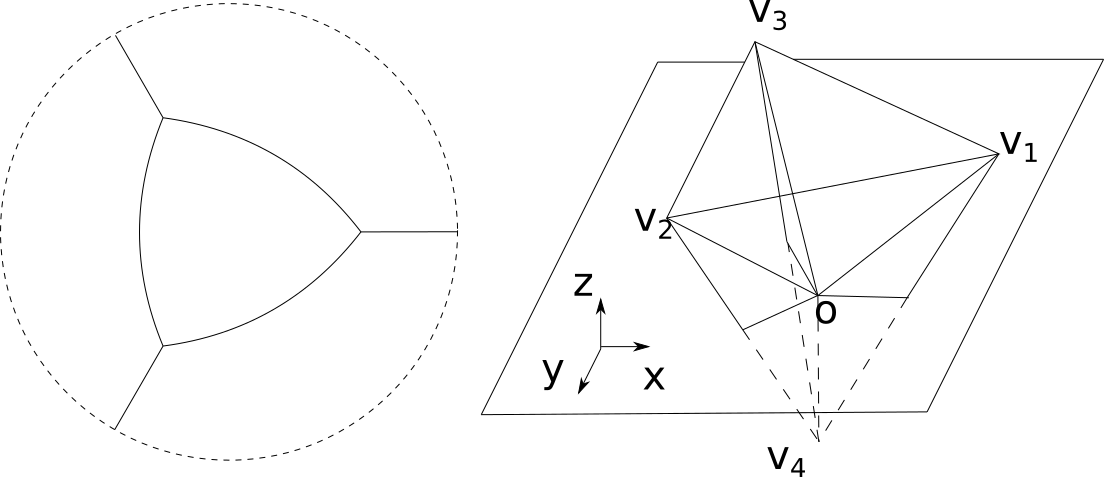}
\caption{On the right the cone $\mathbf{T}_+$ or ``half $\mathbf{T}$'', and on the left its intersection with the hemisphere.}
\label{calibrazioneT}
\end{figure}

The first minimal cone we are going to discuss is called $\T_+$ (or half $\T$) and is obtained by taking a cone of type $\T$ as in Figure \ref{YeTintro.png}, flipping it upside down, placing its barycentre at the origin, and finally intersecting it with the half-space $\R^3_+$ (see Figure \ref{calibrazioneT}).
 In this section we are going to prove the following theorem:

\begin{thm}\label{teoremaT+}
The cone $\T_+$ is an $\alpha$-sliding minimiser in the half-space $\R^3_+$ with respect to $\Gamma=\partial\R^3_+$ if and only if $\alpha\ge\sqrt{\frac{2}{3}}$.
\end{thm}

Let us begin with a formal description of $\T_+$. Given the following unitary vectors
\begin{equation}
\begin{array}{lrrrl}
v_1=\Big( & 2\frac{\sqrt{2}}{3}, &0, &\frac{1}{3} &\Big)\\
v_2=\Big( & -\frac{\sqrt{2}}{3}, &\sqrt{\frac{2}{3}}, &\frac{1}{3} &\Big)\\
v_3=\Big( & -\frac{\sqrt{2}}{3}, &-\sqrt{\frac{2}{3}}, &\frac{1}{3} &\Big)\\
v_4=\Big( & 0, &0, &-1 &\Big)\\
\end{array}
\end{equation}
let $\Delta^3:=[v_1,v_2,v_3,v_4]$ be the 3-dimensional simplex whose vertices are $v_1,v_2,v_3,v_4$, and let $\mathbf{\Delta}^3:=cone(\sk_{1}(\Delta^3))$; then $\T_+:=\mathbf{\Delta}^3\cap\R^3_+$.

We set $\Delta^3_+:=\Delta^3\cap\R^3_+$; for $1=1,2,3,4$ let $F_i:=\Delta^3_i$  be the two-dimensional face of $\Delta^3$ opposed to the vertex $v_i$  and $F_i^+:=F_i\cap\R^3_+$. Let $M\subset\R^3_+$ be a sliding competitor for $\T_+$ such that the symmetric difference between the two is contained in $\Delta^3_+\cap\textrm{Int}(\Delta^3)$. It follows that $\R^3\setminus M$ has 2 unbounded connected components: one of them contains $F_4^+$, and the other one contains the other three faces $F_1^+$, $F_2^+$, $F_3^+$, as well as the lower half-space $\R^3\setminus\R^3_+$. However $\R^3_+\setminus M$ has 4 unbounded connected components each one of them containing one of the faces $F_i^+$.

For $i=1,2,3,4$ we name $V_i$ the connected component of $\R^3_+\setminus M$ containing $F_i^+$ and we set $V_0:=\Delta^3\setminus\R^3_+$. In case $\R^3_+\setminus M$ also has some bounded connected components we just include them in $V_1$. By the definition of sliding competitor we have that $M$ is a Lipschitz image of $\T_+$, therefore $M$ has locally finite $2$-dimensional Hausdorff measure. This means that any of the $V_i$ is a set whose perimeter is locally finite. In particular the sets $U_i:=V_i\cap\Delta^3$ have finite perimeter (in the following we will apply the divergence theorem on these sets). Let us now introduce the following notation:
\begin{eqnarray}
M_{ij} &:=& \partial^* U_i\cap\partial^* U_j\\
M_i &:=& \bigcup_{j=0,\, j\neq i}^4M_{ij}\\
M_i^+ &:=& \bigcup_{j=1}^4M_{ij}, \textrm{ for } i=1,2,3,4\\
\widetilde{M} &:=& M_1^+\cup M_2^+\cup M_3^+\cup M_4
\end{eqnarray}
where $\partial^*E$ denotes the reduced boundary of $E$ (see \cite{ambrosio2000functions} for a definition) and $i,j=0,1,2,3,4$ unless otherwise specified.  The sets $M_{ij}$ are contained in $M$ and in particular $M_{ij}$ is contained in the interface between the regions $U_i$ and $U_j$. Let us now remark that for every $i=0,...,n$, $\H^{n-1}$-almost every point of $M_i$ lies on the interface between exactly two regions $U_i$ and $U_j$. Therefore the interfaces between different couples of regions are essentially disjoint with respect to $\H^{n-1}$ and
\begin{equation}
\H^{n-1}(M_i)=\H^{n-1}\left(\bigcup_{j\neq i}M_{ij}\right)=\sum_{j\neq i}\H^{n-1}\left(M_{ij}\right).
\end{equation}
In order to show it let us first define, for $i=0,...,n$, the exceptional set
\begin{equation}
E_i:=\R^n\setminus\left(U_i^0\cup\partial^*U_i\cup U_i^1\right),
\end{equation}
where, for an $\mathcal{L}^n$-measurable set $A\subset\R^n$ and $t\in[0,1]$, we define
\begin{equation}
A^t:=\left\{x\in\R^n:\lim_{r\to0}\frac{\mathcal{L}^n(A\cap B_r(x))}{\mathcal{L}^n( B_r(x))}=t\right\}.
\end{equation}
Since $U_i$ is a finite perimeter set for every $i=0,...,n$, by Federer's theorem $\H^{n-1}(E_i)=0$ (see \cite[Theorem 3.61]{ambrosio2000functions}). Therefore the exceptional set $E:=\cup_iE_i$ is negligible with respect to $\H^{n-1}$. Let us now assume that a point $x$ belongs to the common boundaries of at least three sets $U_i$, $U_j$ and $U_k$. Clearly the point $x$ cannot belong to the reduced boundaries of the three of them because in this case the blow-up limit of each one of them in the point $x$ would be a half-space and that is a contradiction. Let us assume $x\notin\partial^*U_i$, then $x\in E_i\subset E$. Since $E$ is $\H^{n-1}$-negligible, the same holds true for its intersection with $M_i$ for any $i=0,...,n$. 

We just showed that $\widetilde{M}\subset M\cap\Delta^3_+$ and $\H^2$-almost every point in $\widetilde{M}$ lies on the interface between exactly two regions of its complement (taking into account also $U_0$). Moreover the interfaces between different couples of regions are essentially disjoint with respect to $\H^2$.

 Let us now remark the following useful facts
\begin{eqnarray}
M_i &=& M_i^+\cup M_{i0}\\
\widetilde{M}\setminus\Gamma &=& \bigcup_{i,j\neq0}M_{ij}\quad \textrm{ (up to $\H^2$-negligible sets)}\\
\widetilde{M}\cap\Gamma &=& M_{40}\quad \textrm{ (up to $\H^2$-negligible sets)}\\
\partial^*U_i &=& F_i^+\cup M_i=F_i^+\cup M_i^+\cup M_{i0}\\
\partial^*U_0 &=& \partial^*(\Delta^3\setminus\R^3_+)=(\partial^*\Delta^3\setminus\R^3_+)\cup M_0.
\end{eqnarray}
Finally we denote with $n_i$ the exterior unit normal to $\partial U_i$, and $n_{ij}$ will denote the unit normal to $M_{ij}$ pointing in direction of $U_j$. The vectors of the calibration we will use are $w_i:=\frac{-v_i}{|v_i-v_j|}=-\sqrt{\frac{3}{8}}v_i$ for $i=1,2,3,4$; whose components are the following:
\begin{equation}\label{calibrazione_mezzoT}
\begin{array}{lrrrl}
w_1=\Big( & -\frac{1}{\sqrt{3}}, &0, &-\frac{1}{2\sqrt{6}} &\Big)\\
w_2=\Big( & \frac{1}{2\sqrt{3}}, &-\frac{1}{2}, &-\frac{1}{2\sqrt{6}} &\Big)\\
w_3=\Big( & \frac{1}{2\sqrt{3}}, &\frac{1}{2}, &-\frac{1}{2\sqrt{6}} &\Big)\\
w_4=\Big( & 0, &0, &\frac{1}{2}\sqrt{\frac{3}{2}} &\Big).\\
\end{array}
\end{equation}
Let us remark that

We are now ready to start the paired calibration machinery. After applying the divergence theorem to the sets $U_i$ with the constant vectorfields $w_i$ we can isolate the interface with the negative half-space as follows
\begin{equation}\label{calcoli_calibrazione_mezzoT}
\begin{aligned}
\sqrt{\frac{3}{8}}\H^2\left(\cup_iF^+_i\right)
&=\sum_{i=1}^4\int_{F_i^+}w_i\cdot n_id\mathcal{H}^2\\
&=-\sum_{i=1}^4\int_{M_i}w_i\cdot n_id\mathcal{H}^2\\
&=-\sum_{i=1}^4\int_{M_i^+}w_i\cdot n_id\mathcal{H}^2-\sum_{i=1}^4\int_{M_{i0}}w_i\cdot n_{i0}d\mathcal{H}^2\\
&=\sum_{1\le i<j\le4}\int_{M_{ij}}(w_j-w_i)\cdot n_{ij}d\mathcal{H}^2+\sum_{i=1}^4\int_{M_{i0}}w_i\cdot \hat{z}d\mathcal{H}^2,
\end{aligned}
\end{equation}
where we denoted with $\hat{z}$ the vector $(0,0,1)$. We could replace $n_{i0}$ with $-\hat{z}$ because all the interfaces of kind $M_{i0}$ are contained in the horizontal plane $\Gamma$ and their exterior normal points downward. Let us now focus on the second of the two sums. Using \eqref{calibrazione_mezzoT} and the definition of $M_{ij}$ we get
\begin{equation}
\begin{aligned}\label{calcoli_interfaccia}
\sum_{i=1}^4\int_{M_{i0}}w_i\cdot \hat{z}d\mathcal{H}^2 &=-\sum_{i=1}^3\frac{1}{2\sqrt{6}}\mathcal{H}^2(M_{i0})+\sqrt{\frac{3}{8}}\mathcal{H}^2(M_{40})\\
&=-\frac{1}{2\sqrt{6}}\mathcal{H}^2(M_{10}\cup M_{20}\cup M_{30})+\sqrt{\frac{3}{8}}\mathcal{H}^2(M_{40}).
\end{aligned}
\end{equation}
Plugging \eqref{calcoli_interfaccia} in \eqref{calcoli_calibrazione_mezzoT} and using the fact that $\mathcal{H}^2(M_{10}\cup M_{20}\cup M_{30})=\mathcal{H}^2(M_0)-\mathcal{H}^2(M_{40})$ we obtain

\begin{equation}
\begin{aligned}
&\sqrt{\frac{3}{8}}\mathcal{H}^2(\cup_iF_i^+)+\frac{1}{2\sqrt{6}}\mathcal{H}^2(M_0)=\\
&=\sum_{1\le i<j \le4}\int_{M_{ij}}(w_j-w_i)\cdot n_{ij}d\mathcal{H}^2+\sqrt{\frac{2}{3}}\mathcal{H}^2(M_{40})\\
&\le\H^2\left(\widetilde{M}\setminus\Gamma\right)+\sqrt{\frac{2}{3}}\H^2\left(\widetilde{M}\cap\Gamma\right).
\end{aligned}
\end{equation}
The previous inequality follows from the fact that, by definition $|w_i-w_j|=1$ hence $|(w_j-w_i)\cdot n_{ij}|\le1$.

Let us now assume $\alpha=\sqrt{\frac{2}{3}}$. Recalling that $\widetilde{M}\subset M\cap\Delta^3$ the previous inequality becomes
\begin{equation}\label{calibrazione_quasi_fatto}
\sqrt{\frac{3}{8}}\mathcal{H}^2(\cup_iF_i^+)+\frac{1}{2\sqrt{6}}\mathcal{H}^2(M_0)\le J_\alpha(\widetilde{M})\le J_\alpha(M).
\end{equation}
Since the left-hand side of \eqref{calibrazione_quasi_fatto} is a constant, we just found a lower bound for the energy of any sliding competitor to the cone. Let us now remark that, if we choose $M=\T_+$, we have that $w_i-w_j=n_{ij}$, implying $(w_j-w_i)\cdot n_{ij}\equiv1$. Therefore the inequality \ref{calibrazione_quasi_fatto} turns into an equality in this case. Hence what we just proved is that
\begin{displaymath}
J_\alpha(\T_+)\le J_\alpha(M),
\end{displaymath}
when $\alpha=\sqrt{\frac{2}{3}}$ and for any $M$ sliding competitor to $\T_+$.

This argument also entails that the cone is minimal for $\alpha'\ge\sqrt{\frac{2}{3}}$. It is due to the fact that $J_{\alpha'}(\T_+)=J_{\alpha}(\T_+)$ because $\H^2(\T_+\cap\Gamma)=0$. To show the $\alpha'$-minimality of $\T_+$ we can compute as follows:
\begin{equation}
J_{\alpha'}(\T_+)=J_{\alpha}(\T_+)\le J_{\alpha}(M)\le J_{\alpha'}(M)
\end{equation}
for every sliding competitor $M$.

\begin{figure}
\centering
\includegraphics[scale=0.45]{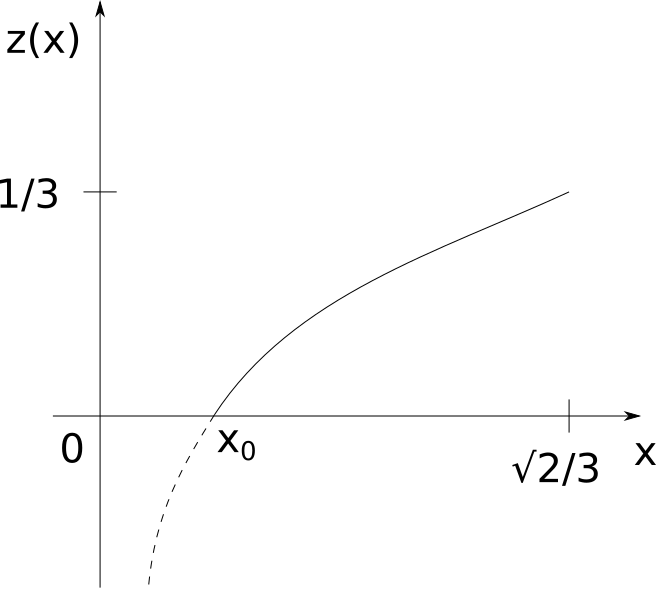}
\caption{Graph of the profile function $z$, the dotted line is its negative part.}
\label{funzionez}
\end{figure}

Let us now check that $\mathbf{T}_+$ is not a sliding minimiser for every $\alpha<\sqrt{\frac{2}{3}}$, we will do it by providing a better competitor. The competitor we are about to show is similar to the competitor for the cone over the skeleton of a cube  provided by Brakke in \cite{brakke1991minimal} and it is a modification of $\mathbf{T}_+$ obtained by pushing it down on $\Gamma$ in such a way to create a little horizontal equilateral triangle centred at the origin and bending the sloping folds along a profile defined by the positive part of the following function (see Figure \ref{funzionez})

\begin{figure}
\centering
\includegraphics[scale=0.5]{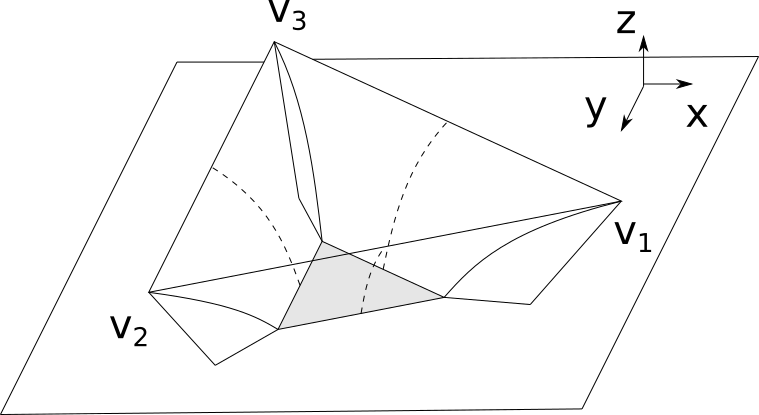}
\caption{Competitor $M$. The dotted line represent the profile given by the function \eqref{profilo_competitore}. The grey region is the intersection between $M$ and $\Gamma$.}
\label{competitoreT}
\end{figure}

\begin{equation}\label{profilo_competitore}
z(x)=\frac{x}{\sqrt{2}}+c\log\left(\frac{3}{\sqrt{2}}x\right),\quad z'(x)=\frac{1}{\sqrt{2}}+\frac{c}{x}.
\end{equation}
The function has been chosen in such a way that $z(\sqrt{2}/3)=1/3$ for any $c>0$ (which will be fixed later). Moreover we set $x_0\in(0,\sqrt{2}/3)$ as the unique solution of $z(x_0)=0$
\begin{equation}\label{eqx0}
0=\frac{x_0}{\sqrt{2}}+c\log\left(\frac{3}{\sqrt{2}}x_0\right).
\end{equation}
What we just defined is a one-parameter family of admissible competitors for $\T_+$. It is possible to use either $c$ or $x_0$ as parameter of the family, and in particular we can use \eqref{eqx0} in order to write $c$ in terms of $x_0$ as follows
\begin{equation}\label{cx0}
c=-\frac{x_0}{\sqrt{2}\log\left(\frac{3}{\sqrt{2}}x_0\right)}.
\end{equation}

Let us now show that for every $\alpha<\sqrt{\frac{2}{3}}$ there exist a competitor in this family, denoted by $M_c$, with less energy than the cone. As usual we define $M_c^+:=M_c\setminus\Gamma$, and we have

\begin{equation}\label{JalphaMc}
J_\alpha(M_c)=\H^2(M_c^+)+\alpha\H^2(M\cap \Gamma).
\end{equation}
By construction the part of $M_c$ laying on $\Gamma$ is just an equilateral triangle whose apothem is $x_0$ therefore its area is $3\sqrt{3}x_0^2$. On the other hand $M_c^+$ is composed by three equal vertical folds and three equal curved folds. Let us call $B$ and $V$ respectively any of the bended or vertical folds, hence
\begin{equation}
\H^2(M_c^+)=3\H^2(B)+3\H^2(V),
\end{equation}
and \eqref{JalphaMc} becomes
\begin{equation}\label{JalphaMc2}
J_\alpha(M_c)=3\H^2(B)+3\H^2(V)+\alpha\H^2(M\cap \Gamma).
\end{equation}
Since the profile of both $B$ and $V$ can be described in terms of the function $z(x)$, their area can be computed by slicing along the direction of the $x$ axis and then integrating on the interval $[x_0,\sqrt{2}/3]$ (or, up to a symmetry, on the interval $[-x_0,-\sqrt{2}/3]$). We have that
\begin{equation}
\H^2(B)=\int_{x_0}^{\frac{\sqrt{2}}{3}}2\sqrt{3}x\sqrt{1+(z'(x))^2}dx,
\end{equation}
where $2\sqrt{3}x$ is the length of the slice and $\sqrt{1+(z'(x))^2}$ is the Jacobian of the function $Z(x):=(x,z(x))$; and
\begin{equation}
\H^2(V)=\int_{x_0}^{\frac{\sqrt{2}}{3}}2z(x)dx
\end{equation}
where we have to multiply by two because the fold $V$ that we are slicing makes and angle $\pi/3$ with the direction of the $x$ axis.
Using the inequality
\begin{equation}
\sqrt{\frac{3}{2}x^2+\sqrt{2}cx+c^2}\le\sqrt{\frac{3}{2}}x+\frac{c}{\sqrt{3}}+\sqrt{\frac{2}{3}}\frac{c^2}{x}
\end{equation}
we can now compute as follows
\begin{equation}\label{faldaB}
\begin{aligned}
\H^2(B) &=2\sqrt{3}\int_{x_0}^{\frac{\sqrt{2}}{3}}x\sqrt{1+\frac{1}{2}+\sqrt{2}\frac{c}{x}+\frac{c^2}{x^2}}dx\\
&=2\sqrt{3}\int_{x_0}^{\frac{\sqrt{2}}{3}}\sqrt{\frac{3}{2}x^2+\sqrt{2}cx+c^2}dx\\
&\le2\sqrt{3}\int_{x_0}^{\frac{\sqrt{2}}{3}}\left[\sqrt{\frac{3}{2}}x+\frac{c}{\sqrt{3}}+\sqrt{\frac{2}{3}}\frac{c^2}{x}\right]dx\\
&=2\sqrt{3}\left[\frac{1}{2}\sqrt{\frac{3}{2}}x^2+\frac{c}{\sqrt{3}}x+\sqrt{\frac{2}{3}}c^2\log\left(\frac{3}{\sqrt{2}}x\right)\right]_{x_0}^\frac{\sqrt{2}}{3}
\end{aligned}
\end{equation}
and
\begin{equation}\label{faldaV}
\begin{aligned}
\H^2(V) &=2\int_{x_0}^{\frac{\sqrt{2}}{3}}\left[\frac{x}{\sqrt{2}}+c\log\left(\frac{3}{\sqrt{2}}x\right)\right]dx\\
&=2\left[\frac{x^2}{2\sqrt{2}}+cx\log\left(\frac{3}{\sqrt{2}}x\right)-cx\right]_{x_0}^\frac{\sqrt{2}}{3}.
\end{aligned}
\end{equation}
Therefore, using \eqref{JalphaMc2}, \eqref{faldaB} and \eqref{faldaV} we obtain
\begin{equation}
J_\alpha(M_c)\le6\left[\sqrt{2}x^2+\sqrt{2}c^2\log\left(\frac{3}{\sqrt{2}}x\right)+cx\log\left(\frac{3}{\sqrt{2}}x\right)\right]_{x_0}^\frac{\sqrt{2}}{3}+\alpha3\sqrt{3}x_0^2.
\end{equation}
Let us now compute the energy of the cone $\mathbf{T}_+$. Since $\mathbf{T}_+=M_0$ we can just use the previous computation with $c=0$ and $x_0=0$ and we get
\begin{equation}
\begin{aligned}
z(x) &=\frac{x}{\sqrt{2}},\quad z'(x)=\frac{1}{\sqrt{2}},\\
J_\alpha(\mathbf{T}_+) &=3\int_{0}^{\frac{\sqrt{2}}{3}}2\sqrt{3}x\sqrt{1+(z'(x))^2}dx+3\int_{0}^{\frac{\sqrt{2}}{3}}2z(x)dx\\
&=12\sqrt{2}\int_{0}^{\frac{\sqrt{2}}{3}}xdx=\frac{4}{3}\sqrt{2}.
\end{aligned}
\end{equation}
Now we can compare the energy of the competitor with the energy of the cone and, using  \eqref{cx0}, we have
\begin{equation}
\begin{aligned}
&J_\alpha(M_c)-J_\alpha(M_0) =\\
=&6\left[\sqrt{2}\frac{2}{9}-\!\sqrt{2}x_0^2-\!\sqrt{2}c^2\log\left(\frac{3}{\sqrt{2}}x_0\right)-\!cx_0\log\left(\frac{3}{\sqrt{2}}x_0\right)\right]+\!\alpha3\sqrt{3}x_0^2-\!\frac{4}{3}\sqrt{2}\\
=&-6\sqrt{2}x_0^2-6\sqrt{2}c^2\log\left(\frac{3}{\sqrt{2}}x_0\right)-6cx_0\log\left(\frac{3}{\sqrt{2}}x_0\right)+\alpha3\sqrt{3}x_0^2\\
=&-6\sqrt{2}x_0^2-3\sqrt{2}\frac{x_0^2}{\log\left(\frac{3}{\sqrt{2}}x_0\right)}+3\sqrt{2}x^2_0+\alpha3\sqrt{3}x_0^2\\
=&3x_0^2\left[-\sqrt{2}-\frac{\sqrt{2}}{\log\left(\frac{3}{\sqrt{2}}x_0\right)}+\alpha\sqrt{3}\right].
\end{aligned}
\end{equation}
Therefore the competitor has less energy than the cone if
\begin{equation}
\alpha\le\sqrt{\frac{2}{3}}\left[1+\frac{1}{\log\left(\frac{3}{\sqrt{2}}x_0\right)}\right].
\end{equation}
Since
\begin{equation}
\lim_{x_0\to0^+}\left[1+\frac{1}{\log\left(\frac{3}{\sqrt{2}}x_0\right)}\right]=1^-
\end{equation}
it follows that for every $\alpha<\sqrt{\frac{2}{3}}$ there exists an $x_0$ (and hence a $c$) such that $J_\alpha(M_c)\le J_\alpha(\mathbf{T}_+)$. This completes the proof of \ref{teoremaT+}.

\begin{figure}
\centering
\includegraphics[scale=0.5]{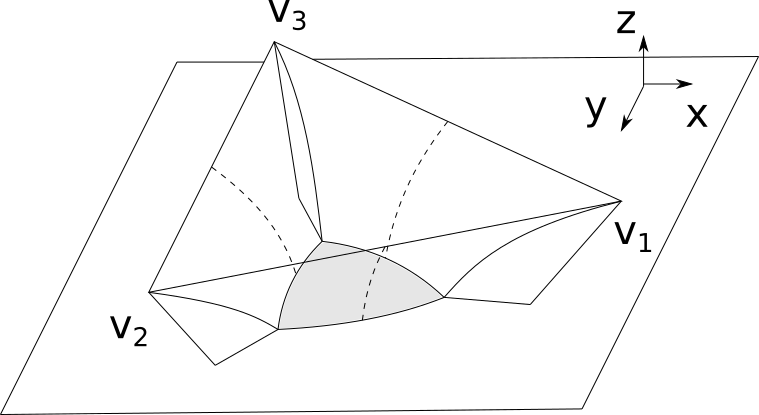}
\caption{Minimiser for $\alpha<\overline{\alpha}$ obtained with Brakke's Surface Evolver. The grey region (fat triangle) is the intersection between the set and $\Gamma$.}
\label{minimoT}
\end{figure}

Let us now set $\overline{\alpha}=\sqrt{\frac{2}{3}}$. In order to understand why $\overline{\alpha}$ is the threshold between minimality and non-minimality of the cone let us recover how it shows up in the calibration argument as the difference of two scalar products
\begin{equation}\label{valore_soglia}
\overline{\alpha}=\sqrt{\frac{2}{3}}=\frac{1}{2}\sqrt{\frac{3}{2}}+\frac{1}{2\sqrt{6}}=w_4\cdot\hat{z}-w_i\cdot\hat{z}.
\end{equation}
The last term in \eqref{valore_soglia} can be reformulated as follows
\begin{equation}
w_4\cdot\hat{z}-w_i\cdot\hat{z}=(w_4-w_i)\cdot\hat{z}=n_{i4}\cdot\hat{z},
\end{equation}
where $n_{i4}$ is the unit normal to one of the sloping folds of $\T_+$. Hence $\overline{\alpha}$ turns out to be the cosine of the angle between the two unit vectors $n_{i4}$ and $\hat{z}$, which is the same as the cosine of the angle between the plane containing the interface $M_{i4}$ and $\Gamma$ since the previous vectors are the unit normals to these planes. Therefore when $\alpha=\overline{\alpha}$ the sloping folds of $\T_+$ satisfy the optimal profile condition $\cos\theta_{\overline{\alpha}}=\overline{\alpha}$, stated in Section \ref{Half-plane}. In case $\alpha<\overline{\alpha}$ the corresponding optimal profile angle $\theta_\alpha$ is bigger than $\theta_{\overline{\alpha}}$ and a competitor, in order to minimise its energy, would try to attain such optimal angle with the sloping folds; resulting in a shape similar to $M_c$. Numerical simulations with Brakke's Surface Evolver \cite{brakke2013surface} show that in this case the minimiser is very similar to $M_c$, it particular the part of it laying on $\Gamma$ is a ``fat'' triangle (see Figure \ref{minimoT}).

On the other hand if $\alpha>\overline{\alpha}$ the optimal profile angle $\theta_\alpha$ is smaller than $\theta_{\overline{\alpha}}$ and for a competitor is impossible to attain it with its sloping folds minimising the energy at the same time.

\subsection{$\Y_\beta$}

\begin{figure}
\centering
\includegraphics[scale=0.35]{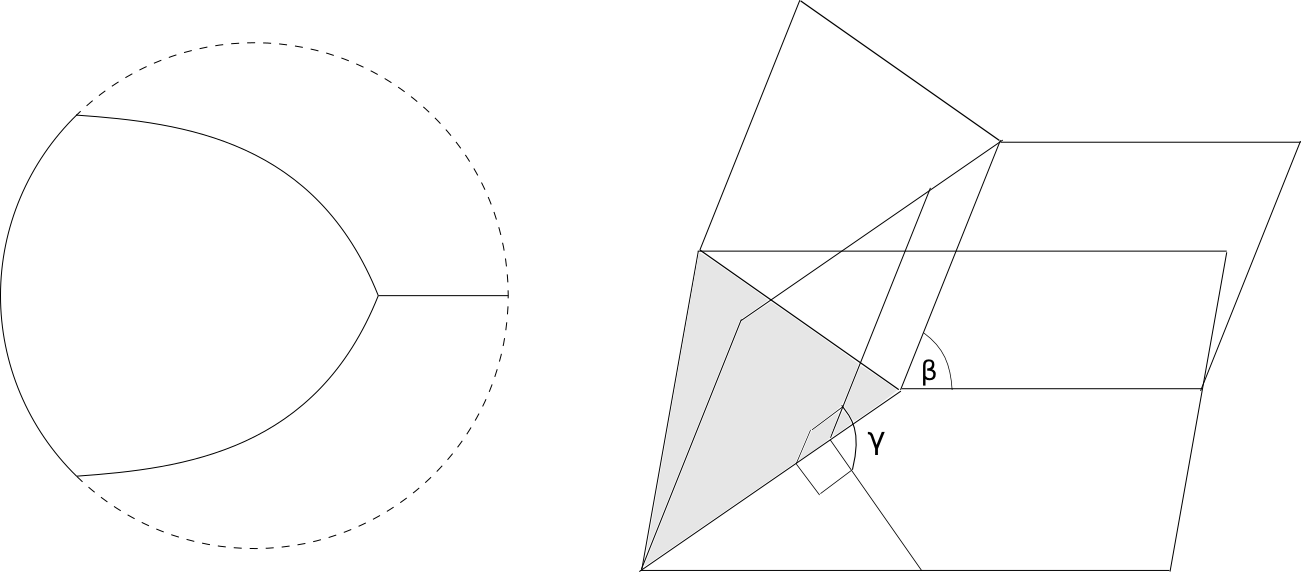}
\caption{The cone $\Y_\beta$ (the grey region is the intersection between the cone and $\Gamma$) and on the left its intersection with the hemisphere.}
\label{Y_tagliato.png}
\end{figure}

Let us introduce a new kind of cone that we will call $\Y_\beta$, where $\beta\in[0,\pi/2]$. It can be obtained with the following procedure: first take the cone $\Y\subset\mathbb{R}^3$, tilt it in a proper way, then intersect it with $\R^3_+$, and finally join it with a section of the horizontal plane $\Gamma$ (see Figure \ref{Y_tagliato.png}).

\begin{figure}
\centering
\includegraphics[scale=0.35]{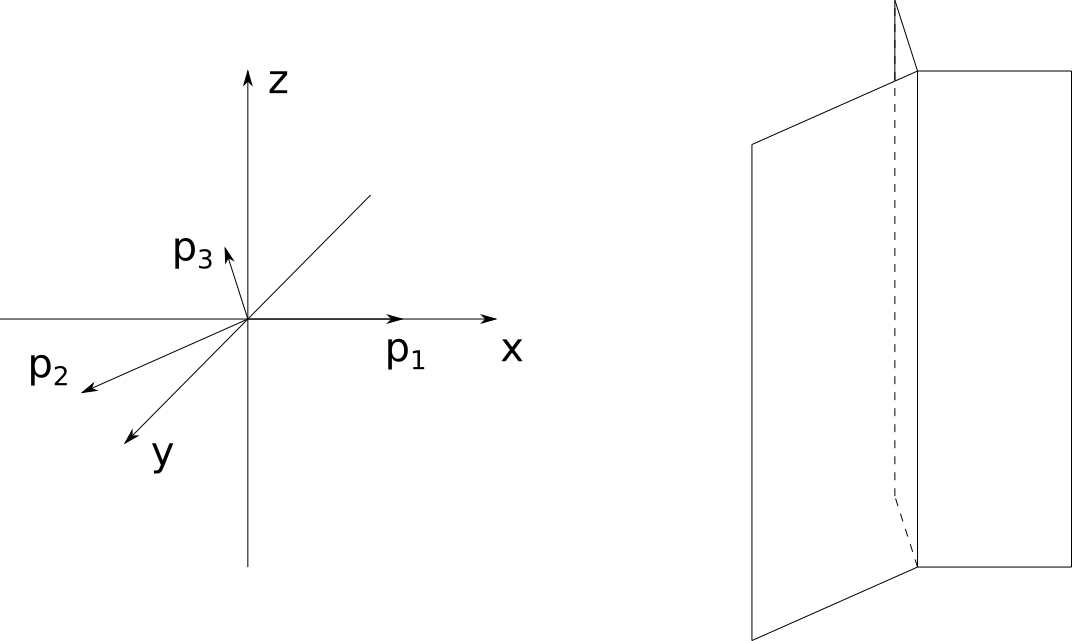}
\caption{The three points generating the cone $Y$ on the left, and the cone $\Y$ on the right.}
\label{yeY.png}
\end{figure}

A formal construction of the cone $\Y_\beta$ is the following. Let $r$ be the straight line spanned by the vector $\hat{z}=(0,0,1)$ and $Y\subset\R^3$ be the 1-dimensional cone over the three points
\begin{equation}
\begin{aligned}
p_1 &=(1,0,0)\\
p_2 &=\left(-\frac{1}{2},\frac{\sqrt{3}}{2},0\right)\\
p_3 &=\left(-\frac{1}{2},-\frac{\sqrt{3}}{2},0\right).
\end{aligned}
\end{equation}
We define the 2-dimensional cone $\Y$ as the Cartesian product $Y\times r$ (see Figure \ref{yeY.png}).
Now we can rotate $\Y$ around the $y$ axis using the following rotation
\begin{equation}
R_\beta:=\left(\begin{array}{ccc}
\sin\beta &0 &\cos\beta\\
0 &1 &0\\
-\cos\beta &0 &\sin\beta
\end{array}\right)
\end{equation}
and we obtain a cone, $R_\beta(\Y)$, composed of one vertical fold and two sloping ones meeting $\Gamma$ with the same angle $\gamma$ (by symmetry). Let us now introduce the vectors
\begin{equation}
\begin{aligned}
n_2 &:=\left(-\frac{\sqrt{3}}{2},-\frac{1}{2},0\right)\\
n_3 &:=\left(-\frac{\sqrt{3}}{2},\frac{1}{2},0\right),
\end{aligned}
\end{equation}
which are respectively orthogonal to $p_2$ and $p_3$ and are contained in $\Gamma$. It follows that the normal vectors to the sloping folds of $R_\beta(\Y)$ can be obtained by rotating $n_2$ and $n_3$ with $R_\beta$; that is to say:
\begin{equation}
\begin{aligned}
m_2:=R_\beta(n_2)=&\left(-\frac{\sqrt{3}}{2}\sin\beta,-\frac{1}{2},\frac{\sqrt{3}}{2}\cos\beta\right)\\
m_3:=R_\beta(n_3)=&\left(-\frac{\sqrt{3}}{2}\sin\beta,\frac{1}{2},\frac{\sqrt{3}}{2}\cos\beta\right).
\end{aligned}
\end{equation}
Since $\cos\gamma=(m_2,\hat{z})$ we find the following relation between $\beta$ and $\gamma$:
\begin{equation}\label{condizione_Yb}
\frac{\sqrt{3}}{2}\cos\beta=\cos\gamma.
\end{equation}
The intersection between $R_\beta(\Y)$ and $\Gamma$ is the union of three half-lines meeting at the origin, each one being the intersection of one of the folds with $\Gamma$. We name $q_1$, $q_2$ and $q_3$ these half-lines, and using $m_2$ and $m_3$ we find that
\begin{equation}\label{retteq123}
\begin{aligned}
q_1 &=\{(t,0,0): t\ge0\}\\
q_2 &=\left\{\left(t,-\sqrt{3}\sin\beta \;t,0\right): t\le0\right\}\\
q_3 &=\left\{\left(t,\sqrt{3}\sin\beta \;t,0\right): t\le0\right\}.
\end{aligned}
\end{equation}
Let us call $S$ the convex subset of $\Gamma$ bounded by $q_2\cup q_3$; we can now define our cone as
\begin{equation}
\Y_\beta:=(R_\beta(\Y)\cap\R^3_+)\cup S.
\end{equation}

In the the rest of this section we will prove the following theorem.
\begin{thm}\label{teoremaYbeta}
The cone $\Y_\beta$ is sliding minimal if and only if
\begin{equation}\label{condizioneYb}
\alpha=\frac{\sqrt{3}}{2}\cos\beta.
\end{equation}
\end{thm}
In order for $\Y_\beta$ to be a minimal set a necessary condition is that for every $P\in\Y_\beta$ the blow-up of $\Y_\beta$ at $P$ has to be a minimal cone (here we are assuming $P\neq0$ because otherwise the necessary condition would turn into a tautology since the blow-up of $\Y_\beta$ at the origin is $\Y_\beta$ itself). In particular if $P\in q_1\cup q_2\cup q_3$ the blow-up of $\Y_\beta$ has to assume one of the optimal profiles described in Section \ref{Half-plane}. In case $P\in q_1$ this condition is satisfied because the cone assumes the profile (ii). When $P\in q_2\cup q_3$ the profile taken by the cone is of type (iv), and in this case the blow-up is a minimal cone if and only if
\begin{equation}
\cos\gamma=\alpha.
\end{equation}
Hence, combining the equality \eqref{condizione_Yb} with the previous one we obtain that the condition \eqref{condizioneYb} expressed in Theorem \ref{teoremaYbeta} is necessary for the minimality of $\Y_\beta$.

We are now going to prove with a calibration argument that condition \eqref{condizioneYb} is also sufficient for the minimality of $\Y_\beta$. The calibration we will use is obtained by rotating with $R_\beta$ a calibration for the cone $Y\subset\Gamma$. Let
\begin{equation}
\begin{aligned}
v_1 &:=\left(-\frac{1}{\sqrt{3}},0,0\right)\\
v_2 &:=\left(\frac{1}{2\sqrt{3}},-\frac{1}{2},0\right)\\
v_3 &:=\left(\frac{1}{2\sqrt{3}},\frac{1}{2},0\right)\\
\end{aligned}
\end{equation}
be a calibration for $Y$ in $\Gamma$ (hence for $\Y$ in $\R^3$), we define the calibration for $\Y_\beta$ as $w_i:=R_\beta(v_i)$ and we get
\begin{equation}
\begin{aligned}
w_1 &=\left(-\frac{\sin\beta}{\sqrt{3}},0,\frac{\cos\beta}{\sqrt{3}}\right)\\
w_2 &=\left(\frac{\sin\beta}{2\sqrt{3}},-\frac{1}{2},-\frac{\cos\beta}{2\sqrt{3}}\right)\\
w_3 &=\left(\frac{\sin\beta}{2\sqrt{3}},\frac{1}{2},-\frac{\cos\beta}{2\sqrt{3}}\right).
\end{aligned}
\end{equation}

\begin{figure}
\centering
\includegraphics[scale=0.35]{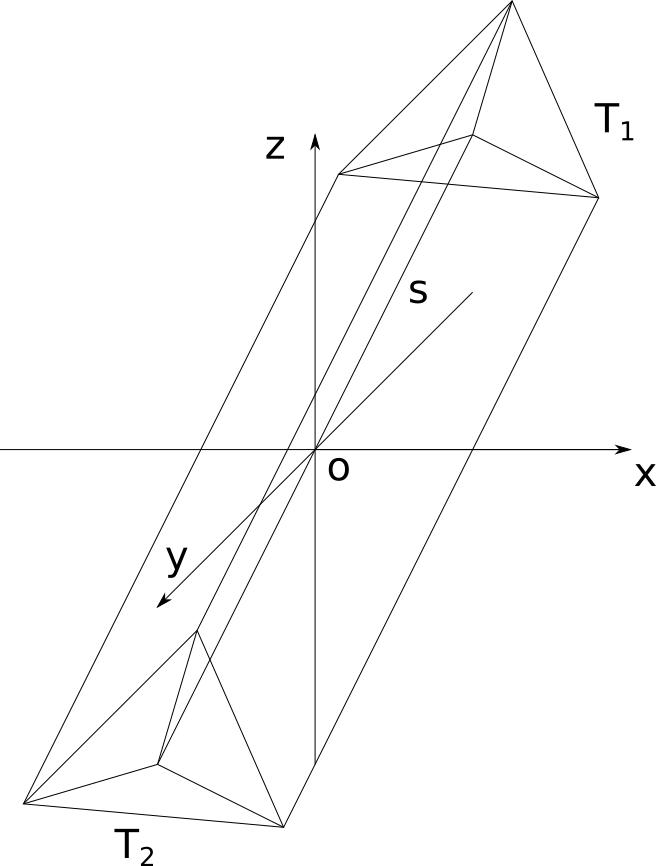}
\caption{The right prism $P$ enclosing part of the cone $R_\beta(\Y)$.}
\label{prisma.png}
\end{figure}

Let $s$ be the straight line spanned by $R_\beta(\hat{z})=(\cos\beta,0,\sin\beta)$, in the following we will refer to it as the spine of $\Y_\beta$. We fix the compact set in which the sliding deformation takes place as the right prism $P$ whose bases are two equilateral triangles, $T_1$ and $T_2$, orthogonal to the spine  $s$ and centred on it, such that their vertices lie in $\Y_\beta$. We assume that the barycentre of $P$ is the origin and its height is large enough such that the two triangular bases do not intersect $\Gamma$ (see Figure \ref{prisma.png}).

By definition the vectors $w_i$ are orthogonal to the lateral faces of the prism $P$, therefore, for $i=1,2,3$ we name $F_i$ the face orthogonal to $w_i$. We set $P_+:=P\cap\R^3_+$ and $F_i^+:=F_i\cap\R^3_+$. Let $M\subset\R^3_+$ be a sliding competitor for $\Y_\beta$ such that $M\triangle\Y_\beta\subset P_+$. It follows that $\R^3_+\setminus M$ has 3 unbounded connected components. For $i=1,2,3$ we name $V_i$ the connected component of $\R^3_+\setminus M$ containing $F_i^+$, and we set $V_0:=\Delta^3\setminus\R^3_+$. In case $\R^3_+\setminus M$ also has bounded connected components we can just include them in $V_1$. For $i=0,1,2,3$ the sets $V_i$ have locally finite perimeter, hence the sets $U_i:=V_i\cap P$ are finite perimeter sets. We can now introduce the following notation
\begin{eqnarray}
M_{ij} &:=& \partial^* U_i\cap\partial^* U_j\\
M_i &:=& \bigcup_{j=0, \, j\neq i}^3M_{ij}\\
M_i^+ &:=& \bigcup_{j=1}^3M_{ij}, \textrm{ for } 1=1,2,3\\
\widetilde{M} &:=& M_1^+\cup M_2^+\cup M_3.
\end{eqnarray}
As in the previous subsection it follows that $\widetilde{M}\subset M\cap P$, $\H^2$-almost every point in $\widetilde{M}$ lies on the interface between exactly two of the $U_i$, and the sets $M_{ij}$ are essentially disjoint when $i<j$. Finally we call $n_i$ the outer normal to $\partial U_i$, and $n_{ij}$ the unit normal to $M_{ij}$ pointing in direction of $U_j$.

We are now ready for the calibration argument. In particular, when applying the divergence theorem to the sets $U_i$ with respect to the vectors $w_i$, we can ignore the contribution given by the upper base of the prism since by definition the vectors $w_i$ are orthogonal to the normal vector to $T_1$. Thus we can compute as follows:

\begin{equation}\label{calcoli_calibrazioneYb}
\begin{aligned}
\frac{1}{\sqrt{3}}\sum_{i=1}^3\H^2(F_i^+) =& \sum_{i=1}^3\int_{F_i^+}w_i\cdot n_id\mathcal{H}^2\\
=&-\sum_{i=1}^3\int_{M_i}w_i\cdot n_id\mathcal{H}^2\\
=&-\sum_{i=1}^3\int_{M_i^+}w_i\cdot n_id\mathcal{H}^2-\sum_{i=1}^3\int_{M_{i0}}w_i\cdot n_id\mathcal{H}^2\\
=&\sum_{1\le i<j\le3}\int_{M_{ij}}(w_j-w_i)\cdot n_{ij}d\mathcal{H}^2+\sum_{i=1}^3\int_{M_{i0}}w_i\cdot \hat{z}d\mathcal{H}^2.
\end{aligned}
\end{equation}
In the last line we isolated the contribution given by the interface with the negative half-space, let us now focus on this term. Using the definition of the vectors $w_i$ we get

\begin{equation}\label{calcoli_interfacciaYb}
\sum_{i=1}^3\int_{M_{i0}}w_i\cdot \hat{z}d\mathcal{H}^2 =\frac{\cos\beta}{\sqrt{3}}\mathcal{H}^2(M_{10})-\frac{\cos\beta}{2\sqrt{3}}\left(\mathcal{H}^2(M_{20})+\mathcal{H}^2(M_{30})\right).
\end{equation}
Plugging \eqref{calcoli_interfacciaYb} in \eqref{calcoli_calibrazioneYb} and using the fact that $\mathcal{H}^2(M_{20})+\mathcal{H}^2(M_{30})=\mathcal{H}^2(M_0)-\mathcal{H}^2(M_{10})$ we obtain
\begin{equation}
\begin{aligned}
\frac{1}{\sqrt{3}}\mathcal{H}^2(\cup_iF_i^+)+\frac{\cos\beta}{2\sqrt{3}}\mathcal{H}^2(M_0) &=\sum_{1\le i<j \le3}\int_{M_{ij}}(w_j-w_i)\cdot n_{ij}d\mathcal{H}^2\\
&+\frac{\sqrt{3}}{2}\cos\beta\mathcal{H}^2(M_{10})\\
&\le\H^2\left(\widetilde{M}\setminus\Gamma\right)+\frac{\sqrt{3}}{2}\cos\beta\left(\widetilde{M}\cap\Gamma\right).
\end{aligned}
\end{equation}
Assuming $\alpha=\frac{\sqrt{3}}{2}\cos\beta$ we have
\begin{equation}\label{calibrazione_quasi_fattoYb}
C(\alpha,P)\le J_\alpha(\widetilde{M})\le J_\alpha(M)
\end{equation}
where $C(\alpha,K)$ is a constant that depends only on the parameter $\alpha$ and on the compact set $K$ containing $M\triangle \Y_\beta$, in our case $K=P$.
Since the left-hand side of \eqref{calibrazione_quasi_fattoYb} is a constant, and the chain of inequalities turns into a chain of equalities for $M=\Y_\beta$, we proved that this cone is minimal when $\alpha=\frac{\sqrt{3}}{2}\cos\beta$. Therefore condition \eqref{condizioneYb} is both necessary and sufficient for the cone $\Y_\beta$ to be minimal.

I the same way it has happened with the cone $\T_+$, it might look surprising how a necessary condition for minimality turns into a sufficient one. However, as we did before, let us remark how the constant $\frac{\sqrt{3}}{2}\cos\beta$ showed up from the computation. It appears as the scalar product between the normal vector to the sloping folds and the normal vector to the domain of the sliding boundary. Since the two vectors have unitary norm their scalar product simply is the cosine of the angle between them, which is the same as the cosine of the angle between the planes they are orthogonal to,
\begin{equation}
\frac{\sqrt{3}}{2}\cos\beta=\frac{\cos\beta}{2\sqrt{3}}+\frac{\cos\beta}{\sqrt{3}}=(w_i-w_3,\hat{z})=\cos\gamma.
\end{equation}
And this means that the reason why we impose \eqref{condizioneYb} as necessary for the minimality of the cone, is actually the same reason that makes it sufficient (the optimal profile angle between the sloping folds and $\Gamma$).

 Let us also remark that the cones of type $\Y_\beta$ satisfying that condition form a one parameter family of minimal cones, depending on the angle $\beta\in[0,\pi/2]$ or, equivalently, on the parameter $\alpha\in[0,1]$. In particular when $\alpha=0$ we have $\beta=\pi/2$, therefore $\Y_\beta$ becomes the union of a vertical half $\Y$ with the section of $\Gamma$ contained in between the two half-lines (cfr. \ref{retteq123})
\begin{equation}
\begin{aligned}
q_2 &=\left\{\left(t,-\sqrt{3}\sin\beta \;t,0\right): t\le0\right\}\\
q_3 &=\left\{\left(t,\sqrt{3}\sin\beta \;t,0\right): t\le0\right\}.
\end{aligned}
\end{equation}
However, since in this case the energy functional $J_0$ does not take into account any set laying on $\Gamma$, up to a $J_0$-negligible set $\Y_{\pi/2}$ is the same as a vertical half $\Y$. On the opposite, when $\alpha=1$ we have that $\beta=0$ and $\Y_\beta$ turns into $\mathbf{V}_{\pi/6}:=V_{\pi/6}\times\R$.

\subsection{$\overline{\Y}_\beta$}

\begin{figure}
\centering
\includegraphics[scale=0.4]{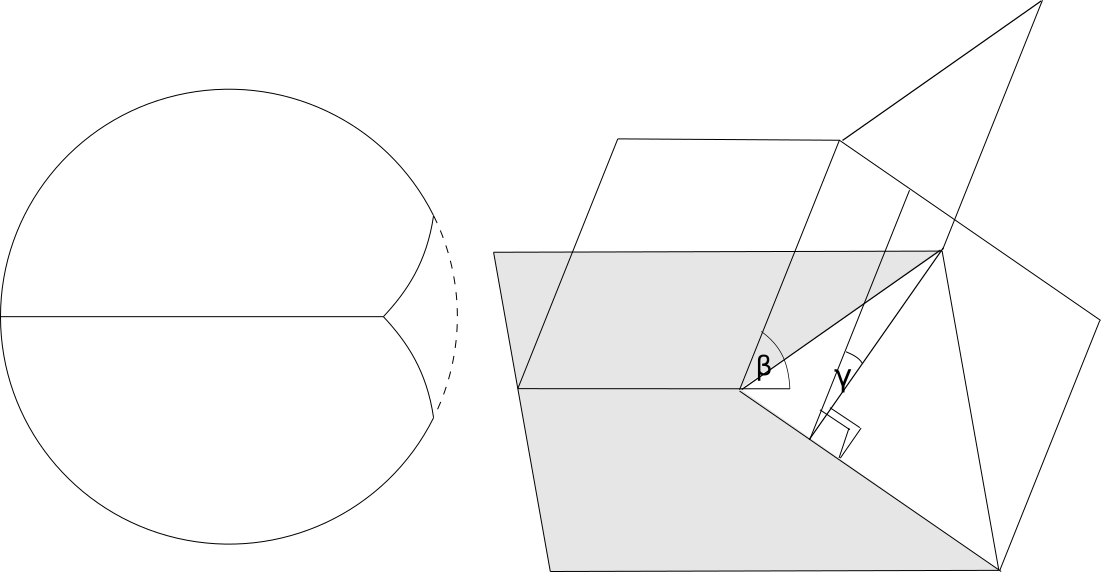}
\caption{The cone $\overline{\Y}_\beta$ (the grey region is the intersection between the cone and $\Gamma$) and on the left its intersection with the hemisphere.}
\label{altro_Y}
\end{figure}

We can use the previous construction to produce another cone that we call $\overline{\Y}_\beta$, for $\beta\in[0,\pi/2]$. It can be obtained with the same procedure as before: first take a cone $\Y\subset\mathbb{R}^3$ symmetric to the previous one, tilt it in a proper way, then intersect it with $\R^3_+$, and finally add to it a section of the horizontal plane $\Gamma$ (see Figure \ref{altro_Y}).

\begin{figure}
\centering
\includegraphics[scale=0.35]{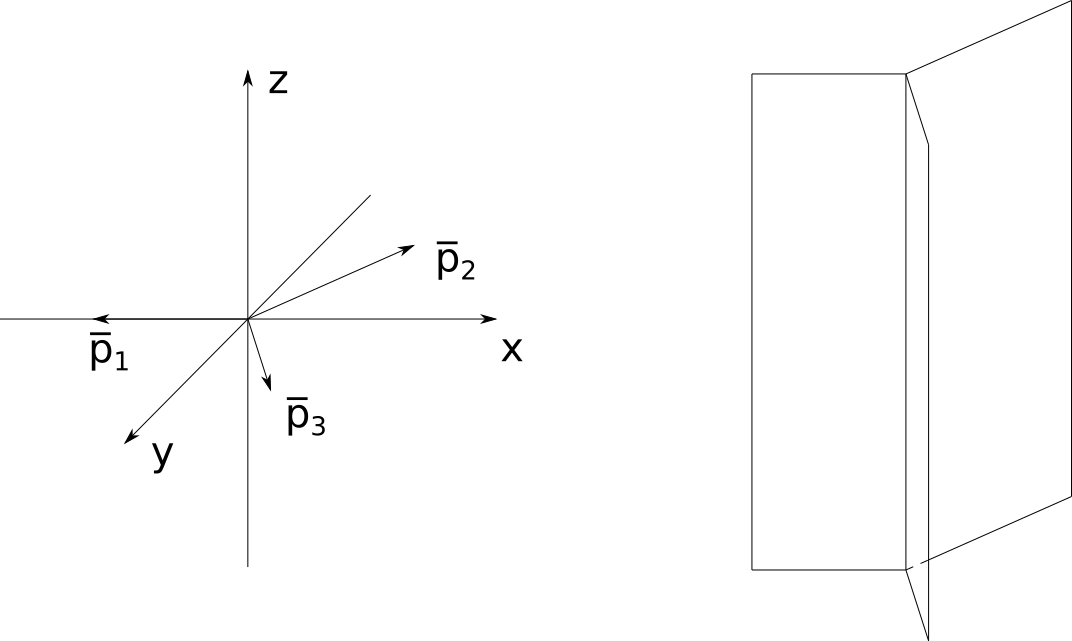}
\caption{The three points generating the cone $\overline{Y}$ on the left, and the cone $\overline{\Y}$ on the right.}
\label{yeY-.png}
\end{figure}

A formal construction of the cone $\overline{\Y}_\beta$ is the following. Let $r$ be the straight line spanned by the vector $\hat{z}=(0,0,1)$ and $\overline{Y}\subset\R^3$ be the 1-dimensional cone over the three points
\begin{equation}
\begin{aligned}
\overline{p}_1 &=(-1,0,0)\\
\overline{p}_2 &=\left(\frac{1}{2},-\frac{\sqrt{3}}{2},0\right)\\
\overline{p}_3 &=\left(\frac{1}{2},\frac{\sqrt{3}}{2},0\right),
\end{aligned}
\end{equation}
in particular $\overline{p}_i=-p_i$ where the $p_i$ are the points defined in the previous subsection, and it follows that $\overline{Y}=-Y$ (where given $E\subset\R^3$ we denote $-E:=\{x\in\R^3:-x\in E\}$). Then we define the 2-dimensional cone $\overline{\Y}$ as the Cartesian product $\overline{Y}\times r$ (see Figure \ref{yeY-.png}), and again we have $\overline{\Y}=-\Y$.
Let us now rotate $\overline{\Y}$ around the $y$ axis with the rotation $R_\beta$. We obtain again a cone, $R_\beta(\overline{\Y})$, composed by one vertical fold and two sloping ones meeting $\Gamma$ with the same angle $\gamma$. The normal vectors to $\overline{p}_2$ and $\overline{p}_3$, respectively $\overline{n}_2$ and $\overline{n}_3$, after an appropriate rotation provide the normal vectors to the sloping folds of $R_\beta(\overline{\Y})$, respectively $\overline{m}_2$ and $\overline{m}_3$. That is
\begin{equation}
\begin{aligned}
\overline{n}_2 &:=\left(-\frac{\sqrt{3}}{2},-\frac{1}{2},0\right)\\
\overline{n}_3 &:=\left(-\frac{\sqrt{3}}{2},\frac{1}{2},0\right)\\
\overline{m}_2 &:=R_\beta(\overline{n}_2)=\left(-\frac{\sqrt{3}}{2}\sin\beta,-\frac{1}{2},\frac{\sqrt{3}}{2}\cos\beta\right)\\
\overline{m}_3 &:=R_\beta(\overline{n}_3)=\left(-\frac{\sqrt{3}}{2}\sin\beta,\frac{1}{2},\frac{\sqrt{3}}{2}\cos\beta\right).
\end{aligned}
\end{equation}
As before, using the fact that $\cos\gamma=(m_2,\hat{z})$, we find the following relation between $\beta$ and $\gamma$:
\begin{equation}\label{condizione_Yb-}
\frac{\sqrt{3}}{2}\cos\beta=\cos\gamma.
\end{equation}
The intersection between $R_\beta(\Y)$ and $\Gamma$ is the union of three half-lines meeting at the origin, each one being the intersection of one of the three folds with $\Gamma$. We name $\overline{q}_1$, $\overline{q}_2$ and $\overline{q}_3$ these half-lines, and using $\overline{m}_2$ and $\overline{m}_3$ we find that
\begin{equation}
\begin{aligned}
\overline{q}_1 &=\{(t,0,0): t\le0\}\\
\overline{q}_2 &=\left\{\left(t,-\sqrt{3}\sin\beta \;t,0\right): t\ge0\right\}\\
\overline{q}_3 &=\left\{\left(t,\sqrt{3}\sin\beta \;t,0\right): t\ge0\right\}.
\end{aligned}
\end{equation}
Let us call $\overline{S}$ the non convex subset of $\Gamma$ bounded by $\overline{q}_2\cup\overline{q}_3$; we can now define our cone as
\begin{equation}
\overline{\Y}_\beta:=(R_\beta(\overline{\Y})\cap\R^3_+)\cup \overline{S}.
\end{equation}

In the the rest of this subsection we will prove the following theorem
\begin{thm}\label{teoremaYbeta-}
The cone $\overline{\Y}_\beta$ is sliding minimal if and only if
\begin{equation}\label{condizioneYb-}
\alpha=\frac{\sqrt{3}}{2}\cos\beta.
\end{equation}
\end{thm}
First of all we have to check that the blow-up of $\overline{\Y}_\beta$ at any of its point $P$ (except for the origin) is a minimal cone. In case $P\in\overline{q}_1$ this condition is satisfied because the cone assumes the profile (iii). When $P\in \overline{q}_2\cup \overline{q}_3$ the profile taken by the cone is of type (iv), and in this case the blow-up is a minimal cone if and only if
\begin{equation}
\cos\gamma=\alpha,
\end{equation}
and we obtain that condition \eqref{condizioneYb-} is necessary for the minimality of $\overline{\Y}_\beta$.

Let us now provide a calibration for $\overline{\Y}_\beta$, this will show that the condition \eqref{condizioneYb-} is also sufficient for the minimality of $\overline{\Y}_\beta$. We will obtain a calibration for $\overline{\Y}_\beta$ by rotating with $R_\beta$ a calibration for the cone $\overline{\Y}$. Let
\begin{equation}
\begin{aligned}
v_1 &:=\left(\frac{1}{\sqrt{3}},0,0\right)\\
v_2 &:=\left(-\frac{1}{2\sqrt{3}},\frac{1}{2},0\right)\\
v_3 &:=\left(-\frac{1}{2\sqrt{3}},-\frac{1}{2},0\right)\\
\end{aligned}
\end{equation}
be a calibration for $\overline{\Y}$ in $\R^3$, we define the calibration for $\overline{\Y}_\beta$ as $\overline{w}_i:=R_\beta(v_i)$ and we get
\begin{equation}
\begin{aligned}
\overline{w}_1 &=\left(\frac{\sin\beta}{\sqrt{3}},0,-\frac{\cos\beta}{\sqrt{3}}\right)\\
\overline{w}_2 &=\left(-\frac{\sin\beta}{2\sqrt{3}},\frac{1}{2},\frac{\cos\beta}{2\sqrt{3}}\right)\\
\overline{w}_3 &=\left(-\frac{\sin\beta}{2\sqrt{3}},-\frac{1}{2},\frac{\cos\beta}{2\sqrt{3}}\right).
\end{aligned}
\end{equation}

\begin{figure}
\centering
\includegraphics[scale=0.35]{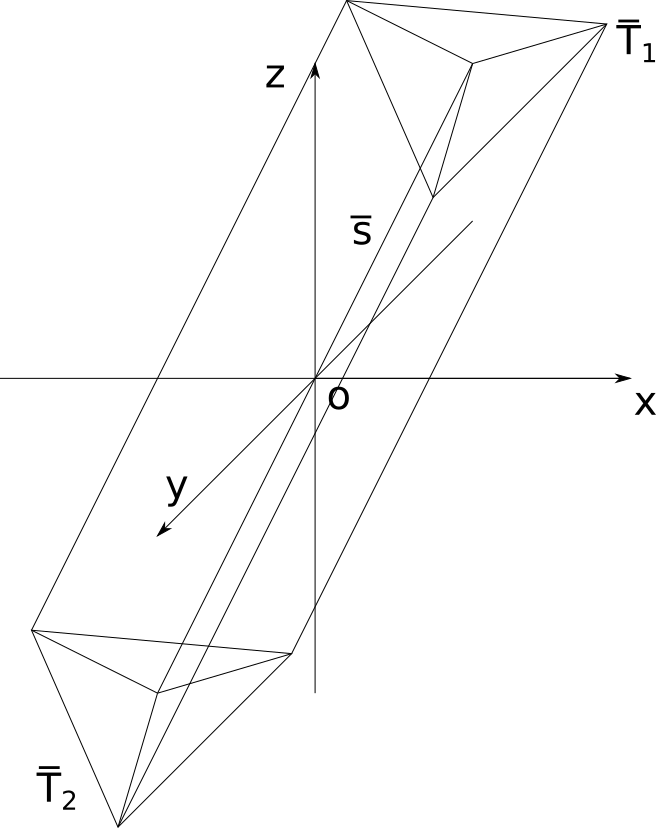}
\caption{The right prism $\overline{P}$ enclosing part of the cone $R_\beta(\overline{\Y})$.}
\label{prisma-.png}
\end{figure}

Let us call $\overline{s}$ the spine of $\overline{\Y}_\beta$, it is spanned by the vector $R_\beta(\hat{z})=(\cos\beta,0,\sin\beta)$. We fix the compact set in which the sliding deformation takes place as the right prism $\overline{P}$ whose bases are two equilateral triangles, $\overline{T}_1$ and $\overline{T}_2$, orthogonal to the spine  $\overline{s}$ and centred on it, such that their vertices lie in $\overline{\Y}_\beta$. We assume that the barycentre of $\overline{P}$ is the origin and its height is large enough such that the two triangular bases do not intersect $\Gamma$ (see Figure \ref{prisma-.png}).

For $i=1,2,3$ we call $F_i$ the lateral face of $\overline{P}$ which is orthogonal to $\overline{w}_i$. We set $\overline{P}_+:=\overline{P}\cap\R^3_+$ and $F_i^+:=F_i\cap\R^3_+$. Let $M\subset\R^3_+$ be a sliding competitor to $\overline{\Y}_\beta$ such that $M\triangle\overline{\Y}_\beta\subset \overline{P}_+$. Than $\R^3_+\setminus M$ has 3 unbounded connected components and we name $V_i$ the one containing $F_i^+$ for $i=1,2,3$, and we set $V_0:=\Delta^3\setminus\R^3_+$. In case $\R^3_+\setminus M$ also has bounded connected components we include them in $V_1$. For $i=0,1,2,3$ the sets $U_i:=V_i\cap \overline{P}$ are finite perimeter sets, and we introduce the following notation:
\begin{eqnarray}
M_{ij} &:=& \partial^* U_i\cap\partial^* U_j\\
M_i &:=& \bigcup_{j=0}^3M_{ij}\\
M_i^+ &:=& \bigcup_{j=1}^3M_{ij}, \textrm{ for } 1=1,2,3\\
\widetilde{M} &:=& M_1^+\cup M_2\cup M_3.
\end{eqnarray}
The set $\widetilde{M}$ defined above is contained in $M$ and $\H^2$-almost every point in it lies on the interface between exactly two of the regions $U_i$; moreover the sets $M_{ij}$ are essentially disjoint. Let us call $n_i$ the outer normal to $\partial U_i$, and $n_{ij}$ the unit normal to $M_{ij}$ pointing in direction of $U_j$. We can now compute as follows

\begin{equation}\label{calcoli_calibrazioneYb-}
\begin{aligned}
\frac{1}{\sqrt{3}}\sum_{i=1}^3\H^2(F_i^+) =& \sum_{i=1}^3\int_{F_i^+}\overline{w}_i\cdot n_id\mathcal{H}^2\\
=&-\sum_{i=1}^3\int_{M_i}\overline{w}_i\cdot n_id\mathcal{H}^2\\
=&-\sum_{i=1}^3\int_{M_i^+}\overline{w}_i\cdot n_id\mathcal{H}^2-\sum_{i=1}^3\int_{M_{i0}}\overline{w}_i\cdot n_id\mathcal{H}^2\\
=&\sum_{1\le i<j\le3}\int_{M_{ij}}(\overline{w}_j-\overline{w}_i)\cdot n_{ij}d\mathcal{H}^2+\sum_{i=1}^3\int_{M_{i0}}\overline{w}_i\cdot \hat{z}d\mathcal{H}^2.
\end{aligned}
\end{equation}
Let us now consider the second term in the last line

\begin{equation}\label{calcoli_interfacciaYb-}
\sum_{i=1}^3\int_{M_{i0}}\overline{w}_i\cdot \hat{z}d\mathcal{H}^2 =-\frac{\cos\beta}{\sqrt{3}}\mathcal{H}^2(M_{10})+\frac{\cos\beta}{2\sqrt{3}}\left(\mathcal{H}^2(M_{20})+\mathcal{H}^2(M_{30})\right).
\end{equation}
The two previous computation together with the fact that $\mathcal{H}^2(M_{10})=\mathcal{H}^2(M_0)-\mathcal{H}^2(M_{20})-\mathcal{H}^2(M_{30})$ lead to
\begin{equation}
\begin{aligned}
\frac{1}{\sqrt{3}}\mathcal{H}^2(\cup_iF_i^+)+\frac{\cos\beta}{\sqrt{3}}\mathcal{H}^2(M_0) &=\sum_{1\le i<j \le3}\int_{M_{ij}}(\overline{w}_j-\overline{w}_i)\cdot n_{ij}d\mathcal{H}^2\\
&+\frac{\sqrt{3}}{2}\cos\beta\left(\mathcal{H}^2(M_{20})+\mathcal{H}^2(M_{30})\right)\\
&\le\H^2\left(\widetilde{M}\setminus\Gamma\right)+\frac{\sqrt{3}}{2}\cos\beta\left(\widetilde{M}\cap\Gamma\right).
\end{aligned}
\end{equation}
Therefore, in case $\alpha=\frac{\sqrt{3}}{2}\cos\beta$, we get
\begin{equation}\label{calibrazione_quasi_fattoYb-}
C(\alpha,P)\le J_\alpha(\widetilde{M})\le J_\alpha(M)
\end{equation}
where $C(\alpha,K)$ is a constant only depending on the parameter $\alpha$ and on the compact set $K$ containing $M\triangle \overline{\Y}_\beta$, in our case $K=P$.
Since the left-hand side of \eqref{calibrazione_quasi_fattoYb-} is a constant, and the chain of inequalities turns into a chain of equalities for $M=\overline{\Y}_\beta$, we proved that this cone is minimal when $\alpha=\frac{\sqrt{3}}{2}\cos\beta$. Therefore condition \eqref{condizioneYb-} is both necessary and sufficient for the cone $\overline{\Y}_\beta$ to be minimal, and once again the explanation of this fact relies on the optimal angle profile.

The minimal cones of type $\overline{\Y}_\beta$ form a one parameter family depending on the angle $\beta\in[0,\pi/2]$ or, equivalently, on the parameter $\alpha\in[0,1]$. In particular when $\alpha=0$ we have $\beta=\pi/2$, therefore $\overline{\Y}_\beta$ becomes the union of a vertical half $\overline{\Y}$ with the section of $\Gamma$ not contained in between the two half-lines
\begin{equation}
\begin{aligned}
\overline{q}_2 &=\left\{\left(t,-\sqrt{3}\sin\beta \;t,0\right): t\ge0\right\}\\
\overline{q}_3 &=\left\{\left(t,\sqrt{3}\sin\beta \;t,0\right): t\ge0\right\}.
\end{aligned}
\end{equation}
However, since in this case the energy functional $J_0$ does not take into account any set laying on $\Gamma$, up to a $J_0$-negligible set $\overline{\Y}_{\pi/2}$ is the same as a vertical half $\overline{\Y}$. On the opposite, when $\alpha=1$ we have that $\beta=0$ and $\overline{\Y}_\beta$ turns into a cone composed by the union of $\Gamma$ with a vertical half-plane.

\subsection{Double $\Y_\beta$}

\begin{figure}
\centering
\includegraphics[scale=0.4]{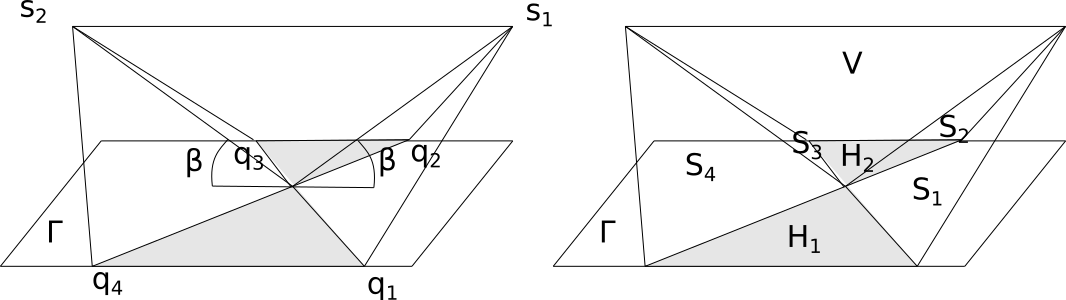}
\caption{The cone $\mathbf{W}_\beta$ with the names of the lines on the left and the names of the folds on the right (the grey region is the intersection between the cone and $\Gamma$).}
\label{doppioY}
\end{figure}

\begin{figure}
\centering
\includegraphics[scale=0.4]{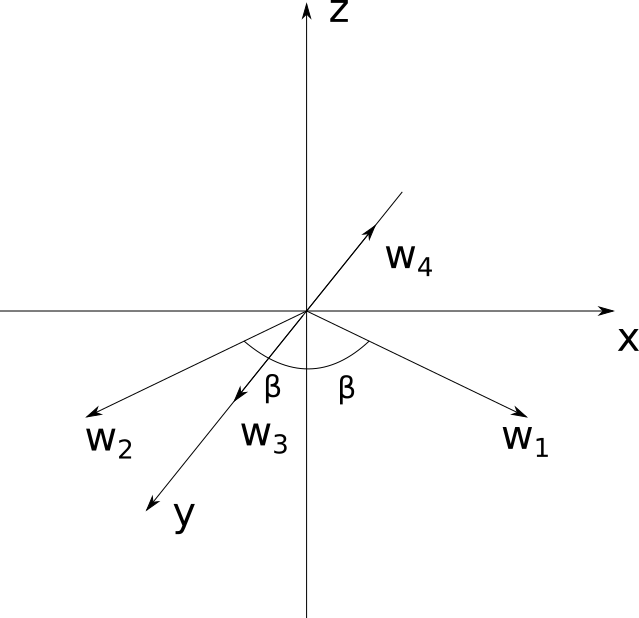}
\caption{Calibration for the cone $\mathbf{W}_\beta$.}
\label{doppioYcal}
\end{figure}

The next cone is called $\mathbf{W}_\beta$, and it is composed by two cones of type $\overline{\Y}_\beta$ symmetric to each other with respect to a vertical plane, and sharing the same vertical fold (see Figure \ref{doppioY}).
\begin{figure}
\centering
\includegraphics[scale=0.35]{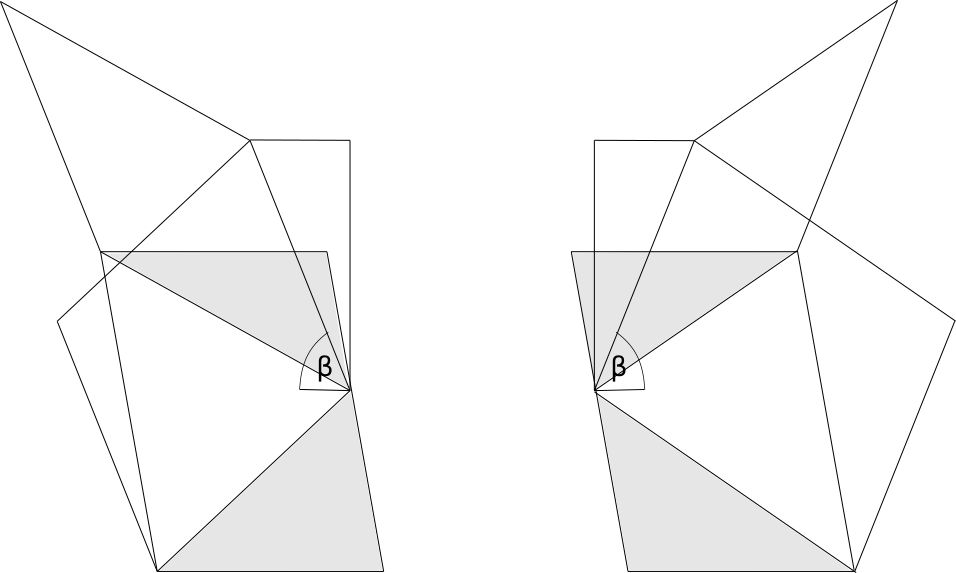}
\caption{The intersection of $\overline{\Y}_\beta$ with $H$ on the right, and its reflection with $R_x$ on the left.}
\label{mezzoYb-}
\end{figure}
It can be constructed as follows. Let $H:=\{(x,y,z)\in\R^3:x\ge0\}$ be the half-space of positive $x$, and $\overline{\Y}_\beta^{x^+}:=\overline{\Y}_\beta\cap H$. Let $R_x$ be the reflection with respect to the $yz$ plane
\begin{equation}
R_x=\left(\begin{array}{ccc}
-1 & 0 & 0\\
0 & 1 & 0\\
0 & 0 & 1\\
\end{array}\right)
\end{equation}
then we can define $W_\beta:=\overline{\Y}_\beta^{x^+}\cup R_x(\overline{\Y}_\beta^{x^+})$ (see Figure \ref{mezzoYb-}).
\begin{figure}
\centering
\includegraphics[scale=0.35]{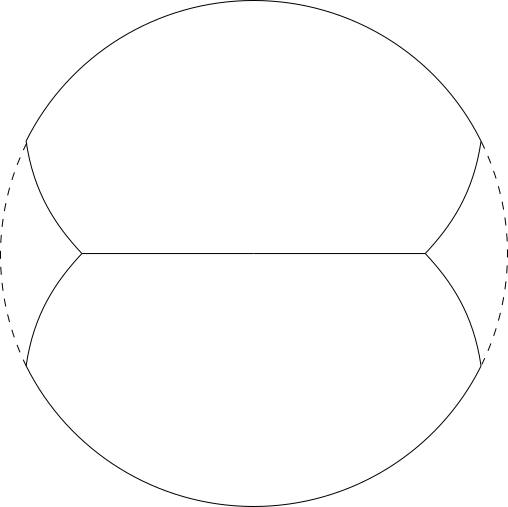}
\caption{The intersection of $\overline{\Y}_\beta$ with the hemisphere.}
\label{intersezione-doppioY}
\end{figure}
In this section we will prove the following theorem.
\begin{thm}\label{teoremaWbeta}
Let $\sin\beta\le1/\sqrt{3}$, then the cone $\mathbf{W}_\beta$ is sliding minimal if and only if
\begin{equation}\label{condizioneWb}
\alpha=\frac{\sqrt{3}}{2}\cos\beta.
\end{equation}
\end{thm}
As usual the necessity of condition \eqref{condizioneWb} is given by the minimality of the tangent cone of $\mathbf{W}_\beta$ at any of its point and the sufficiency is proved via calibration. By construction the two spines of the sloping $\Y$ cones are the two following half lines
\begin{equation}
\begin{aligned}
s_1 &=\{(\cos\beta\;t,0,\sin\beta\;t):t\ge0\}\\
s_2 &=\{(-\cos\beta\; t,0,\sin\beta\; t):t\ge0\},\\
\end{aligned}
\end{equation}
and the intersection of the four sloping folds with $\Gamma$ are the following four half lines
\begin{equation}
\begin{aligned}
q_1 &=\{(t,\sqrt{3}\sin\beta\;t,0):t\ge0\}\\
q_2 &=\{(t,-\sqrt{3}\sin\beta\;t,0):t\ge0\}\\
q_3 &=\{(t,\sqrt{3}\sin\beta\;t,0):t\le0\}\\
q_4 &=\{(t,-\sqrt{3}\sin\beta\;t,0):t\le0\}.
\end{aligned}
\end{equation}
We can name the folds of $\mathbf{W}_\beta$ as follows (see Figure \ref{doppioY}):
\begin{itemize}
\item[$V$:] the vertical planar face bounded by the two spines $s_1$ and $s_2$;
\item[$H_1$:] the horizontal planar face bounded by $q_1$ and $q_4$;
\item[$H_2$:] the horizontal planar face bounded by $q_2$ and $q_3$;
\item[$S_1$:] the sloping planar face bounded by $q_1$ and $s_1$;
\item[$S_2$:] the sloping planar face bounded by $q_2$ and $s_1$;
\item[$S_3$:] the sloping planar face bounded by $q_3$ and $s_2$;
\item[$S_4$:] the sloping planar face bounded by $q_4$ and $s_2$.
\end{itemize}
For $i=1,2,3,4$ let $\hat{s}_i$ be a unit vector orthogonal to $S_i$. Exploiting the computation of the previous section and the symmetries of $\mathbf{W}_\beta$ we get
\begin{equation}
\begin{aligned}
\hat{s}_1 &=\left(-\frac{\sqrt{3}}{2}\sin\beta, \frac{1}{2}, \frac{\sqrt{3}}{2}\cos\beta\right)\\
\hat{s}_2 &=\left(-\frac{\sqrt{3}}{2}\sin\beta, -\frac{1}{2}, \frac{\sqrt{3}}{2}\cos\beta\right)\\
\hat{s}_3 &=\left(\frac{\sqrt{3}}{2}\sin\beta, -\frac{1}{2}, \frac{\sqrt{3}}{2}\cos\beta\right)\\
\hat{s}_4 &=\left(\frac{\sqrt{3}}{2}\sin\beta, \frac{1}{2}, \frac{\sqrt{3}}{2}\cos\beta\right).\\
\end{aligned}
\end{equation}
Therefore we can choose the following vectors as our calibration (see Figure \ref{doppioYcal})
\begin{equation}
\begin{aligned}
w_1 &=\left(\frac{\sqrt{3}}{2}\sin\beta,0,-\frac{\sqrt{3}}{2}\cos\beta\right)\\
w_2 &=\left(-\frac{\sqrt{3}}{2}\sin\beta,0,-\frac{\sqrt{3}}{2}\cos\beta\right)\\
w_3 &=\left(0,\frac{1}{2},0\right)\\
w_4 &=\left(0,-\frac{1}{2},0\right),
\end{aligned}
\end{equation}
and it is easily seen that, except for $w_1-w_2$, the difference between any two vectors of the calibration is the unit normal to some of the folds of $\mathbf{W}_\beta$. In the following we will explain better the role played by this difference in the calibration argument.
Let us name the 4 connected components of $\R^3_+\setminus\W_\beta$ as follows
\begin{itemize}
\item[$V_1$:] the connected component bounded by $S_1$, $S_2$ and $\Gamma$;
\item[$V_2$:] the connected component bounded by $S_3$, $S_4$ and $\Gamma$;
\item[$V_3$:] the connected component bounded by $H_1$, $S_1$, $V$ and $S_4$;
\item[$V_4$:] the connected component bounded by $H_2$, $S_2$, $V$ and $S_3$.
\end{itemize}
Let $B:=B_1(0)$ be the ball with unitary radius centred at the origin. We choose $B$ as the compact set in which the deformation takes place, and we set $F_i:=V_i\cap\partial B$ for $i=1,2,3,4$ and $B_+:=B\cap\R^3_+$. Let $M$ be a sliding competitor to $\W_\beta$ such that $M\triangle\W_\beta\subset B_+$. It follows that $\R^3_+\setminus M$ has 4 unbounded connected components and, with an abuse of notation, we still call them $V_i$ for $i=1,2,3,4$ (in such a way that these connected components correspond to the previous ones when $M=\W_\beta$). We set $V_0:=\Delta^3\setminus\R^3_+$ and, in case $\R^3_+\setminus M$ also has bounded connected components we include them in $V_1$. The sets $V_i$ have locally finite perimeter, hence the sets $U_i:=V_i\cap B$ are finite perimeter sets. Let us introduce the following sets
\begin{eqnarray}
M_{ij} &:=& \partial^* U_i\cap\partial^* U_j\\
M_i &:=& \bigcup_{j=0, \, j\neq i}^4M_{ij} \\
M_i^+ &:=& \bigcup_{j=1}^4M_{ij}, \textrm{ for } 1=1,2,3\\
\widetilde{M} &:=& M_1^+\cup M_2^+\cup M_3\cup M_4.
\end{eqnarray}
It follows that $\widetilde{M}\subset M\cap B$, $\H^2$-almost every point in $\widetilde{M}$ lies on the interface between exactly two of the $U_i$, and the sets $M_{ij}$ are essentially disjoint. Finally we call $n_i$ the outer normal to $\partial U_i$, and $n_{ij}$ the unit normal to $M_{ij}$ pointing in direction of $U_j$. Thus
\begin{equation}\label{calcoli_calibrazioneW}
\begin{aligned}
\sum_{i=1}^4\int_{F_i}w_i\cdot n_i d\H^2=&-\sum_{i=1}^4\int_{M_i}w_i\cdot n_id\mathcal{H}^2\\
=&-\sum_{i=1}^4\int_{M_i^+}w_i\cdot n_id\mathcal{H}^2-\sum_{i=1}^4\int_{M_{i0}}w_i\cdot n_id\mathcal{H}^2\\
=&\sum_{1\le i<j\le4}\int_{M_{ij}}(w_j-w_i)\cdot n_{ij}d\mathcal{H}^2+\sum_{i=1}^4\int_{M_{i0}}w_i\cdot \hat{z}d\mathcal{H}^2.
\end{aligned}
\end{equation}
Let us consider the second term in the last line
\begin{equation}\label{calcoli_interfacciaW}
\begin{aligned}
\sum_{i=1}^4\int_{M_{i0}}w_i\cdot \hat{z}d\mathcal{H}^2 &=-\frac{\sqrt{3}}{2}\cos\beta\left(\H^2(M_{10})\H^2+(M_{20})\right)\\
&=\frac{\sqrt{3}}{2}\cos\beta\left(\H^2(M_{30})\H^2+(M_{40})\right)-\frac{\sqrt{3}}{2}\cos\beta\H^2(M_0).
\end{aligned}
\end{equation}
Putting together the two previous computation we get
\begin{equation}\label{calcoliW}
\begin{aligned}
&\sum_{i=1}^4\int_{F_i}w_i\cdot n_i d\H^2+\frac{\sqrt{3}}{2}\cos\beta\H^2(M_0)=\\
&=\sum_{1\le i<j\le4}\int_{M_{ij}}(w_j-w_i)\cdot n_{ij}d\mathcal{H}^2+\frac{\sqrt{3}}{2}\cos\beta\left(\H^2(M_{30})\H^2+(M_{40})\right)\\
&\le\H^2\left(\widetilde{M}\setminus\Gamma\right)+\frac{\sqrt{3}}{2}\cos\beta\H^2\left(\widetilde{M}\cap\Gamma\right).
\end{aligned}
\end{equation}
In order for the last inequality to be true we have to impose
\begin{equation}\label{condizionecalW}
|w_i-w_j|\le1\quad\forall1\le i<j\le4.
\end{equation}
As we remarked above, in order to satisfy \eqref{condizionecalW} we only need to check that $|w_1-w_2|\le1$, and this condition leads to
\begin{equation}
\sin\beta\le\frac{1}{\sqrt{3}}.
\end{equation}
Since the left-hand side in the first line of \eqref{calcoliW} is a constant depending only on the shape of the compact set chosen and on the calibration (which only depend on $\alpha$) in the following we will denote it with $C(\alpha,B)$. Therefore, assuming  $\alpha=\frac{\sqrt{3}}{2}\cos\beta$ and $\sin\beta\le\frac{1}{\sqrt{3}}$ \eqref{calcoliW} becomes
\begin{equation}
C(\alpha,B)\le J_\alpha(\widetilde{M})\le J_\alpha(M).
\end{equation}
Let us now remark that $\H^2(M_{12})=0$ when $M=\W_\beta$, then in this case the previous inequalities turn into a chain of equalities and we proved the minimality of $\W_\beta$.

\begin{figure}
\centering
\includegraphics[scale=0.39]{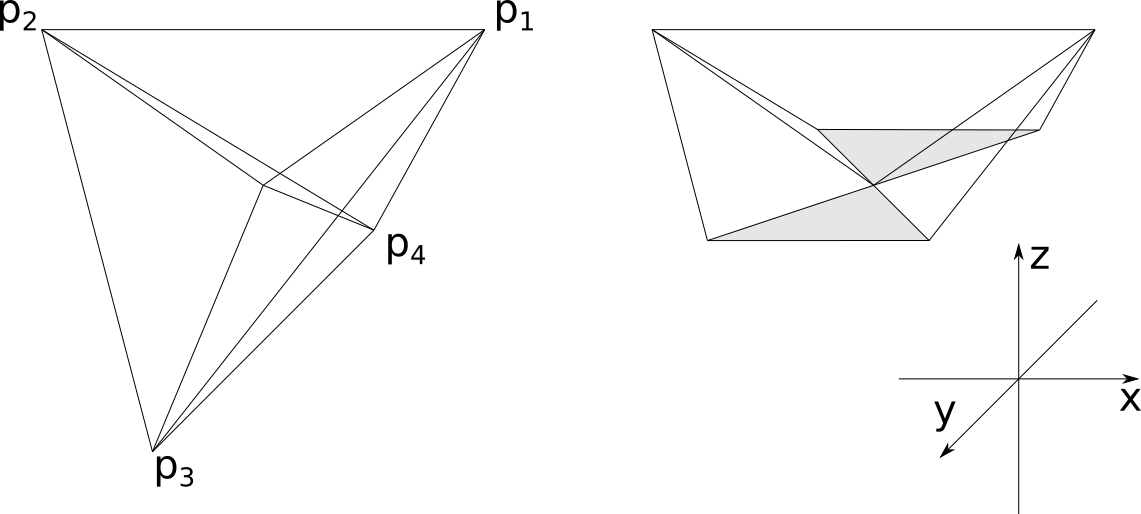}
\caption{On the left the tetrahedron $\Delta_3$ and the cone $\T$ over its skeleton. On the right the corresponding cone $\W_\beta$ and the relative Cartesian coordinate system (not centred at the origin).}
\label{Wtetraedro}
\end{figure}

The cones of type $\W_\beta$ satisfying the minimality condition form a one-parameter family which can be described in therm of the parameter $\alpha\in[1/\sqrt{2},1]$, or equivalently in term of the  angle $\beta \in [0,\arcsin(1/\sqrt{3})]$. In particular, when $\alpha=1$ the cone $\W_\beta$ turns into the union of $\Gamma$ with a vertical half-plane. On the other hand, when $\alpha=1/\sqrt{2}$ the two sloping $\Y$ cones of $\W_\beta$ actually belong to a cone of type $\T$ (see Figure \ref{Wtetraedro}). It can be obtained as the cone over the skeleton of the regular tetrahedron $\Delta_3$ whose vertices are
\begin{equation}
\begin{aligned}
p_1=& \left(\sqrt{\frac{2}{3}},0,\frac{1}{\sqrt{3}}\right) &p_2=& \left(-\sqrt{\frac{2}{3}},0,\frac{1}{\sqrt{3}}\right)\\
p_3=& \left(0,\sqrt{\frac{2}{3}},-\frac{1}{\sqrt{3}}\right) &p_4=& \left(0,-\sqrt{\frac{2}{3}},-\frac{1}{\sqrt{3}}\right).\\
\end{aligned}
\end{equation}

 \bibliography{bibliografia}{}

\begin{thebibliography}{10}

\bibitem{ambrosio2000functions}
Luigi Ambrosio, Nicola Fusco, and Diego Pallara.
\newblock {\em Functions of bounded variation and free discontinuity problems},
  volume 254.
\newblock Clarendon Press Oxford, 2000.

\bibitem{brakke2013surface}
Kenneth~A Brakke.
\newblock Surface {E}volver, {V}ersion 2.70.
\newblock \url{http://facstaff.susqu.edu/brakke/evolver/evolver.html}.
\newblock Accessed: 01-03-2018.

\bibitem{brakke1991minimal}
Kenneth~A Brakke.
\newblock Minimal cones on hypercubes.
\newblock {\em The Journal of Geometric Analysis}, 1(4):329--338, 1991.

\bibitem{david2014local}
Guy David.
\newblock Local regularity properties of almost-and quasiminimal sets with a
  sliding boundary condition.
\newblock {\em arXiv preprint arXiv:1401.1179}, 2014.

\bibitem{fang2016holder}
Yangqin Fang.
\newblock H{\"o}lder regularity at the boundary of two-dimensional sliding
  almost minimal sets.
\newblock {\em Advances in Calculus of Variations}, 2016.

\bibitem{fang2017local}
Yangqin Fang.
\newblock Local ${C}^{1,\beta}$-regularity at the boundary of two dimensional
  sliding almost minimal sets in $\mathbb{R}^3$.
\newblock arXiv:1611.01343, 2017.

\bibitem{federer2014geometric}
Herbert Federer.
\newblock {\em Geometric measure theory}.
\newblock Springer, 2014.

\bibitem{finn1974capillarity}
Robert Finn.
\newblock Capillarity phenomena.
\newblock {\em Russian Mathematical Surveys}, 29(4):133--153, 1974.

\bibitem{giusti1976boundary}
Enrico Giusti.
\newblock Boundary value problems for non-parametric surfaces of prescribed
  mean curvature.
\newblock {\em Annali della Scuola Normale Superiore di Pisa-Classe di
  Scienze}, 3(3):501--548, 1976.

\bibitem{heppes1964isogonal}
A.~Heppes.
\newblock Isogonal sphärischen net.
\newblock {\em Ann. Univ. Sci. Budapest Eötvös Sect. Math.}, 7:41--48, 1964.

\bibitem{lamarle1865stabilite}
Ernest Lamarle.
\newblock Sur la stabilit{\'e} des syst{\`e}mes liquides en lames minces.
\newblock {\em M{\'e}moires de l'Acad{\'e}mie Royale des Sciences, des Lettres
  et des Beaux-Arts de Belgique}, 35:1--104, 1865.

\bibitem{lawlor1994paired}
Gary Lawlor and Frank Morgan.
\newblock Paired calibrations applied to soap films, immiscible fluids, and
  surfaces or networks minimizing other norms.
\newblock {\em Pacific Journal of Mathematics}, 166(1):55--83, 1994.

\bibitem{marchese2016steiner}
Andrea Marchese and Annalisa Massaccesi.
\newblock The steiner tree problem revisited through rectifiable g-currents.
\newblock {\em Advances in Calculus of Variations}, 9(1):19--39, 2016.

\bibitem{massaccesi2014currents}
Annalisa Massaccesi.
\newblock {\em Currents with coefficients in groups, applications and other
  problems in Geometric Measure Theory}.
\newblock PhD thesis, Ph. D. thesis, Scuola Normale Superiore di Pisa, 2014.

\bibitem{dephilippis2015regularity}
G~De Philippis and Francesco Maggi.
\newblock Regularity of free boundaries in anisotropic capillarity problems and
  the validity of young’s law.
\newblock {\em Archive for Rational Mechanics and Analysis}, 216(2):473--568,
  2015.

\bibitem{plateau1873statique}
Joseph Antoine~Ferdinand Plateau.
\newblock {\em Statique exp{\'e}rimentale et th{\'e}orique des liquides soumis
  aux seules forces mol{\'e}culaires}, volume~2.
\newblock Gauthier-Villars, 1873.

\bibitem{taylor1973regularity}
Jean~E Taylor.
\newblock Regularity of the singular sets of two-dimensional area-minimizing
  flat chains modulo 3 inr 3.
\newblock {\em Inventiones mathematicae}, 22(2):119--159, 1973.

\bibitem{taylor1976structure}
Jean~E Taylor.
\newblock The structure of singularities in soap-bubble-like and soap-film-like
  minimal surfaces.
\newblock {\em Annals of Mathematics}, pages 489--539, 1976.

\bibitem{taylor1977boundary}
Jean~E Taylor.
\newblock Boundary regularlty for solutions to various capillarity and free
  boundary problems.
\newblock {\em Communications in Partial Differential Equations},
  2(4):323--357, 1977.

\bibitem{white1999rectifiability}
Brian White.
\newblock Rectifiability of flat chains.
\newblock {\em Annals of Mathematics}, 150:165--184, 1999.

\end{thebibliography}
 \bibliographystyle{plain}
\end{document}